\documentclass[aos,dvips]{arximspdf}
\usepackage{graphics}

\doi{10.1214/09-AOS763}
\volume{38}
\issue{5}
\pubyear{2010}
\firstpage{3129}
\lastpage{3163}

\makeatletter

\newtheorem{theorem}{Theorem}
\newtheorem{corollary}{Corollary}

\newtheorem{lemma}{Lemma}

\newproclaim{assumption}{Assumption}

\makeatother

\begin{document}
\begin{frontmatter}

\title{Nonparametric tests of the Markov hypothesis
in continuous-time models\thanksref{T1}}
\thankstext{T1}{Supported in part by NSF Grants SES-0850533,  DMS-05-32370 and DMS-09-06482.}
\runtitle{Testing the Markov hypothesis}

\begin{aug}
\author[A]{\fnms{Yacine} \snm{A\"{\i}t-Sahalia}\ead[label=e1]{yacine@princeton.edu}},
\author[B]{\fnms{Jianqing} \snm{Fan}\ead[label=e2]{jqfan@princeton.edu}}
\and
\author[C]{\fnms{Jiancheng} \snm{Jiang}\ead[label=e3]{jjiang1@uncc.edu}\corref{}}
\runauthor{Y. A\"{\i}t-Sahalia, J. Fan and J. Jiang}
\affiliation{Princeton University and NBER, Princeton University and University of North~Carolina at Charlotte}
\address[A]{Y. A\"{\i}t-Sahalia\\ Department of Economics\\
Princeton University \\
Princeton, New Jersey 08544\\
USA\\
and\\
 NBER\\
 1050 Massachusetts Ave.\\
 Cambridge, Massachusetts 02138\\ USA\\ \printead{e1}}
\address[B]{J. Fan\\ Department of ORFE\\
Princeton University\\
Princeton, New Jersey 08544\\USA\\\printead{e2}}
\address[C]{J. Jiang\\ Department of Mathematics\\
\quad and Statistics\\
University of North Carolina\\
\quad at Charlotte\\
Charlotte, North Carolina 28223\\
USA\\\printead{e3}} 
\end{aug}

\received{\smonth{6} \syear{2009}}
\revised{\smonth{10} \syear{2009}}

\begin{abstract}
We propose several
statistics to test the Markov hypothesis for $\beta$-mixing stationary
processes sampled at discrete time intervals.
Our tests are based on the Chapman--Kolmogorov equation. We establish
the asymptotic null distributions of the proposed test statistics,
showing that Wilks's phenomenon holds. We compute the power of the test
and provide simulations to investigate
the finite sample performance of the test statistics when the null model is a diffusion
process, with alternatives consisting of models with a stochastic mean reversion level,
stochastic volatility and jumps.
\end{abstract}



\end{frontmatter}

\section{Introduction}\label{sec:intro}

Among stochastic processes, those that satisfy the\break Markov property represent
an important special case. The Markov property restricts the effective size of
the filtration that governs the dynamics of the process. In a nutshell, only
the current value of $X$ is relevant to determine its future evolution. This
restriction simplifies model-building, forecasting and time series inference.
Can it be tested on the basis of discrete observations? It is not practical to
approach the testing problem in the form of a restriction on the
filtration,\ the size of any alternative filtration being essentially
unrestricted. Furthermore, the continuous-time filtration is not observable on
the basis of discrete observations, especially if we do not have high-frequency
data, and asymptotically the sampling interval remains fixed.

Instead, we propose to test the Markov property at the level of\break the
discrete-frequency transition densities of the process. Given a\break
time-homogeneous stochastic process $X=\{X_{t}\}_{t\geq0}$ on $\mathbb{R}^{m}%
$, with the standard probability space $(\Omega;{\mathcal{F}};P)$ and
filtration ${\mathcal{F}}_{t}\subset{\mathcal{F}}$, we consider families of
conditional probability functions $P(\cdot|x,\Delta)$ of $X_{t+\Delta}$ given
$X_{t}=x$: for each Borel measurable function $\psi$, $E[\psi(X_{t+\Delta
})|{\mathcal{F}}_{t}] = \int\psi(y)P(dy|X_{t},\Delta)$.

If $X$ is time-homogeneous Markovian, then its transition densities satisfy
the Chapman--Kolmogorov equation
\begin{equation}\label{2.1}
P(\cdot|x,\Delta+\tau)=\int_{S}P(\cdot|y,\Delta)P(dy|x,\tau)
\end{equation}
for all $\Delta>0$ and $\tau>0$ and $x$ in the support $S$ of $X.$ Suppose
that we collect $n$ observations on $X$ on $[0,T]$ sampled every $\Delta$
units of time. We will assume that $\Delta$ is fixed; asymptotics are
therefore with $T\rightarrow\infty$. High-frequency asymptotics, by contrast,
assume that $\Delta\rightarrow0$, and $T$ can be fixed or $T$ diverges. This
asymptotic setup could have been considered, but it is not necessary here as
we are able to test the hypothesis on the basis of discrete data at a fixed
interval with no requirement for high-frequency data; high-frequency
asymptotics would, of course, also generate different asymptotic properties for
the tests we propose.

If we set $\tau=\Delta$ in (\ref{2.1}), then we can estimate the transition
densities at the desired frequencies on the basis of these discrete
observations. On the left-hand side of the equation, the transition density at
interval $2\Delta$ can be estimated simply by retaining every other
observation in the same data sample. To avoid unnecessary restrictions on the
data-generating process, we will employ nonparametric estimators of the
transition densities. Given these, equation (\ref{2.1}) then becomes a
testable implication of the Markov property for $X$.

Conversely, Kolmogorov's construction (see, e.g., \cite{revuzyor},
Chapter III, Theorem~1.5) allows one to parameterize Markov processes using
transition functions. Namely, given a transition function $P$ and a
probability measure $\pi$ on $\mathbb{R}^{m}$ serving as the initial
distribution, there exists a unique probability measure such that the
coordinate process $X$ is Markovian with respect to $\sigma(X_{u},u\leq t)$,
has transition function $P$ and $X_{0}$ has $\pi$ as its distribution. When
$\pi$ is the invariant probability measure of $P$, the process is a stationary
Markov process. Therefore, given an initial distribution, a Markov process
$X$ is determined by its transition densities.

Transition densities play a crucial role in many contexts. In mathematical
finance, arbitrage considerations in finance make many pricing problems
linear; as a result, they depend upon the computation of conditional
expectations for which knowledge of the transition function is essential.
Also, inference strategies relying on maximum-likelihood or Bayesian methods
require the transition density of the process. Specification testing
procedures for stochastic processes also make use of the transition densities
(see, e.g., \cite{yacrfs96,yacfanpeng09,chengao04b,chengaotong08,gaocasas08} and \cite{hongli05}). All these
models, estimation methods and tests assume that the process is Markovian.

Stochastic volatility models are a very broad class of non-Markovian models,
due to the latency of the volatility state variable. They have been popular in
financial asset pricing and modeling (see, e.g.,
\cite{fouquepapanicolaousircar00}). Parameters in stochastic volatility models
are much harder to estimate and the associated pricing formulas are also
different from those based on Markovian diffusion models and depend on
the assumptions made on the correlation structure between the innovations to
prices and volatility (as in, e.g., \cite{heston93}). Other examples include
models for the term structure of interest rates, which may be Markovian or not
(see, e.g., \cite{heathjarrowmorton92}), and, in fact, one popular approach in
mathematical finance consists of restricting term structure models to be
Markovian (see, e.g., \cite{caverhill94}). In other words, many financial
econometrics models are based on the Markovian assumption and this fundamental
assumption needs to be tested before they can be applied. In all these cases,
testing whether the underlying process is Markovian is essential in helping to
decide which family of models to use and whether a diffusion model is adequate.

We will propose test statistics for this purpose. Asymptotic null
distributions of test statistics are established and we show that Wilks's
phenomenon holds for several of those test statistics. The power functions of
the tests are also computed for contiguous alternatives. We find that the
proposed tests can detect alternatives with an optimal rate in the context of
nonparametric testing procedures.

The remainder of the paper is organized as follows. In Section \ref{sec:non}, we
briefly describe the nonparametric estimation of the transition functions of
the process. In Section \ref{sec:test}, we propose several test statistics for
checking the Markov hypothesis. In Section \ref{sec:asy}, we establish their
asymptotic null distributions and compute their power. Simulation results are
reported in Section \ref{sec:sim}. Technical conditions and proofs of the
mathematical results are given in Section \ref{sec:tech}.

\section{Nonparametric estimation of the transition density and distribution
functions}\label{sec:non}

To estimate nonparametrically the transition density of observed process $X,$
we use the locally linear method suggested by \cite{fanyaotong96}. The process
$X$ is sampled at regular time points $\{i\Delta, i=1,\ldots,n+2\}$. We make
the dependence on the transition function and related quantities on $\Delta$
implicit by redefining
\[
X_{i}=X_{i\Delta},\qquad i=1,\ldots,n+2,
\]
which is assumed to be a stationary and $\beta$-mixing process.

For ease of exposition, we describe the estimation of the transition density
and distribution when $m=1$, that is, $X$ is a process on the line. We also define
$Y_{i}=Y_{i\Delta}=X_{(i+1)\Delta}$ and $Z_{i}=Z_{i\Delta}=X_{(i+2)\Delta}$.
Let $b_{1}$ and $b_{2}$ denote two bandwidths and $K$ and $W$ two kernel
functions. Observe that as $b_{2}\rightarrow0$
\begin{equation} \label{2.9}
E[  K_{b_{2}}(Z_{i}-z)|Y_{i}=y]  \approx p(z|y,\Delta),
\end{equation}
where $K_{b_{2}}(z)=K(z/b_{2})/b_{2}$ and $p(z|y, \Delta)$ is the transition
density of $X_{(i+1)\Delta}$ given $X_{i\Delta}$. The left-hand side of
(\ref{2.9}) is the regression function of the random variable $K_{b_{2}}%
(Z_{i}-z)$ given $Y_{i}=y$. Hence, locally linear fit can be used to estimate
this regression function. For each given $x$, one minimizes
\begin{equation}\label{2.10}
\sum_{i=1}^{n}\{K_{b_{2}}(Z_{i}-z)-\alpha-\beta(Y_{i} -y)\}^{2}W_{b_{1}}%
(Y_{i}-y)
\end{equation}
with respect to the the local parameters $\alpha$ and $\beta$, where
$W_{b_{1}}(z)=W(z/b_{1})/b_{1}$. The resulting estimate of the conditional
density is simply $\hat{\alpha}$. The estimator can be explicitly expressed
as
\begin{equation}\label{2.11}
\hat{p}(z|y,\Delta)=n^{-1}\sum_{i=1}^{n}W_{n}(Y_{i}-y,y;b_{1})K_{b_{2}}%
(Z_{i}-z),
\end{equation}
where $W_{n}$ is the effective kernel induced by the local linear fit.
Explicitly, it is given by
\[
W_{n}(z,y;b_{1})=W_{b_{1}}(z)\frac{s_{n,2}(y)-b_{1}^{-1}zs_{n,1}(y)}%
{s_{n,0}(y)s_{n,2}(y)-s_{n,1}^{2}(y)},
\]
where
\[
s_{n,j}(y)=\frac{1}{n}\sum_{i=1}^{n}\biggl(  \frac{Y_{i}-y}{b_{1}}\biggr)
^{j}W_{b_{1}}(Y_{i}-y).
\]
Note that the effective kernel $W_{n}$ depends on the sampling data points and
the location $y$. This is the key to the design adaptation and location
adaptation property of the locally linear fit.

From (\ref{2.11}), a possible estimate of the transition distribution
$P(z|y,\Delta)=P(Z_{i}<z|Y_{i}=y,\Delta)$ is given by
\[
\hat{P}(z|y,\Delta)=\int_{-\infty}^{z}\hat{p}(t|y,\Delta)\, dt
=\frac{1}{n}
\sum_{i=1}^{n}W_{n}(Y_{i}-y,y;b_{1})\bar{K}\biggl(  \frac{Z_{i} -z}{b_{2}%
}\biggr)  ,
\]
where $\bar{K}(u)=\int_{u}^{\infty}K(t) \,dt$. Let $b_{2}\rightarrow0$, then
\begin{equation}\label{2.12}
\hat{P}(z|y,\Delta)=\frac{1}{n}\sum_{i=1}^{n}W_{n}(Y_{i}-y,y;b_{1}
)I(Z_{i}<z),
\end{equation}
where we drop the term in which $Z_{i}=z$ would contribute the value $\bar
{K}(0)$. This does not affect the asymptotic property of $\hat{P}$. Actually,
(\ref{2.12}) is really the locally linear estimator of the regression
function
\[
P(z|y,\Delta)=E[  I(Z_{i}<z)|Y_{i}=y]  .
\]

\section{Nonparametric tests for the Markov hypothesis in discretely sampled
continuous-time models}\label{sec:test}

The tests we propose are based on the fact that, for $X$ to be Markovian, its
transition function must satisfy the Chapman--Kolmogorov equation in the form
for densities equivalent to (\ref{2.1}),
\begin{equation}\label{2.16}
p(z|x,2\Delta)=r(z|x,2\Delta),
\end{equation}
where%
\begin{equation}\label{2.18}
r(z|x,2\Delta)\equiv\int_{y\in S}p(z|y,\Delta)p(y|x,\Delta)\,dy
\end{equation}
for all $(x,z)\in S^{2}$.

Under time-homogeneity of the process $X$, the Markov hypothesis can then be
tested in the form $H_{0}$ against $H_{1},$ where
\begin{equation} \label{2.17}
\cases{
H_{0}\dvtx p(z|x,2\Delta)-r(z|x,2\Delta)=0 &\quad\mbox{for all} $(x,z)\in
S^{2}$,
\cr
H_{1}\dvtx p(z|x,2\Delta)-r(z|x,2\Delta)\neq0 &\quad\mbox{for
some} $(x,z)\in S^{2}$.
}
\end{equation}
This test corresponds to a nonparametric null hypothesis versus a
nonparametric alternative hypothesis.

Both $p(y|x,\Delta)$ and $p(z|x,2\Delta)$ can be estimated from data sampled
at interval $\Delta$, thanks to time homogeneity. In fact, the successive
pairs of observed data $\{(X_{i},Y_{i})\}_{i=1}^{n+1}$ form a sample from the
distribution with conditional density $p(y|x,\Delta)$ from which the
estimator $\hat{p}(y|x,\Delta)$ can be constructed, and then $\tilde
{r}(z|x,2\Delta)$ as indicated in equation (\ref{2.18}) can be computed.
Meanwhile, the successive pairs $(X_{1},Z_{1}),(X_{2},Z_{2}),\ldots,$ form a
sample from the distribution with conditional density $p(z|x,2\Delta)$ which
can be used to form the direct estimator by drawing a parallel to
(\ref{2.11})
\[
\hat{p}(z|x,2\Delta)=\frac{1}{n}\sum_{i=1}^{n}W_{n}(X_{i}-x,x;h_{1})K_{h_{2}%
}(Z_{i}-z),
\]
where $h_{1}$ and $h_{2}$ are two bandwidths, localizing, respectively, the $x$-
and $z$-domain.

In other words, the test compares a direct estimator of the $2\Delta$-interval
conditional density, $\hat{p} (z|x, 2\Delta)$, to an indirect estimator of the
$2\Delta$-interval conditional density, $\tilde{ r} (z|x, 2\Delta)$, obtained
by (\ref{2.18}). If the process is actually Markovian, then the two estimates
should be close (for some distance measure) in a sense made precise by the use
of the statistical distributions of these estimators.

If, instead of $2\Delta$ transitions, we test the replicability of $j\Delta$
transitions, where $j$ is an integer greater than or equal to $2$, there is no
need to explore all the possible combinations of these $j\Delta$ transitions
in terms of shorter ones $(1,j-1),(1,j-2),\ldots$: verifying equation
(\ref{2.16}), for one combination is sufficient as can be seen by a recursion
argument. In the event of a rejection of $H_{0}$ in (\ref{2.17}), there is no
need to consider transitions of order $j$. In general, a vector of
``transition equalities'' can be tested in a
single pass in a method of moments framework with as many moment conditions as
transition intervals.

We propose two classes of tests for the hypothesis problem (\ref{2.17}) based
on nonparametric estimation of the transition densities and distributions. To
be more specific, since
\begin{equation}\label{2.19}
r(z|x,2\Delta)=E[  p(z|Y_{i},\Delta)|X_{i}=x]  ,
\end{equation}
the function $r(z|x,2\Delta)$ can also be estimated by regressing
nonparametrically $\hat{p}(z|Y_{i},\Delta)$ on $X_{i}$. This avoids
integration in (\ref{2.18}) and makes implementation and theoretical studies
easier. Employing the local linear smoother for (\ref{2.19}), we obtain the
following estimator:
\[
\hat{r}(z|x,2\Delta)=n^{-1}\sum_{i=1}^{n}W_{n}(X_{i}-x,x,h_{3})\hat{p}%
(z|Y_{i},\Delta),
\]
where $h_{3}$ is a bandwidth in this smoothing problem. Under $H_{0}$ in
(\ref{2.17}), the logarithm of the likelihood function is estimated as
\[
\ell(H_{0}) =\sum_{i=1}^{n}\log\hat{r}(Z_{i}|X_{i},2\Delta),
\]
after ignoring the initial stationary density $\pi(X_{1})$. This likelihood
can be compared with
\[
\ell(H_{1})=\sum_{i=1}^{n}\log\hat{p}(Z_{i}|X_{i},2\Delta),
\]
which leads to the generalized likelihood ratio (GLR) test statistic (see
\cite{fanzhangzhang01})
\[
\sum_{i=1}^{n}\log\{\hat{r}(Z_{i}|X_{i},2\Delta)/\hat{p}(Z_{i}|X_{i}%
,2\Delta)\}.
\]
Since the nonparametric regression functions cannot be estimated well when
$(X_{i},Z_{i})$ is in the boundary region, the above GLR test statistic is
reduced to
\[
T_{0}=\sum_{i=1}^{n}\log\{\hat{r}(Z_{i}|X_{i},2\Delta)/\hat{p}(Z_{i}%
|X_{i},2\Delta)\}w^{*}(X_{i},Z_{i}),
\]
where $w^{*}$ is a weight function selected to reduce the influences of the
unreliable estimates in the sparse region. Admittedly, $\ell(H_{1})$ is not
the estimated log-likelihood under $H_{1}$ in (\ref{2.17}), but is used to
create a discrepancy measure. To see this, note that under $H_{0}$, $\hat{p}$
and $\hat{r}$ are approximately the same. By Taylor's expansion, we have
\begin{eqnarray*}
T_{0}  &\approx& \sum_{i=1}^{n} \frac{\hat{p}(Z_{i}|X_{i},2\Delta)-\hat
{r}(Z_{i}|X_{i},2\Delta)} {\hat{p}(Z_{i}|X_{i},2\Delta)}w^{*}(X_{i},Z_{i})
\\
&&{}  +\frac{1}{2}\sum_{i=1}^{n} \biggl\{\frac{\hat{p}(Z_{i}|X_{i},2\Delta
)-\hat{r}(Z_{i}|X_{i},2\Delta)} {\hat{p}(Z_{i}|X_{i},2\Delta)}\biggr\}^{2}%
w^{*}(X_{i},Z_{i}).
\end{eqnarray*}
To avoid unnecessary technicalities, we ignore the first term and consider the
second term
\begin{equation}\label{2.20a}
T_{1}^{\ast}=\sum_{i=1}^{n}\biggl\{  \frac{\hat{p}(Z_{i}|X_{i },2\Delta
)-\hat{r}(Z_{i}|X_{i},2\Delta)}{\hat{p}(Z_{i }|X_{i},2\Delta)}\biggr\}
^{2}w^{*}(X_{i},Z_{i}),
\end{equation}
which is the $\chi^{2}$-type of test statistics. A natural alternative
statistic to $T_{1}^{*}$ is
\begin{equation}\label{2.20}
T_{1}=\sum_{i=1}^{n}\{\hat{p}(Z_{i}|X_{i},2\Delta)-\hat{r}(Z_{i}|X_{i}%
,2\Delta)\}^{2}w(X_{i},Z_{i }).
\end{equation}
The resulting test statistics $T_{1}^{*}$ and $T_{1}$ are discrepancy measures
between $\hat{p}$ and $\hat{r}$ in the $L_{2}$-distance. Discrepancy-measure
based test statistics receive attention and achieve success in the literature.
Other discrepancy norms such as the $L_{\infty}$-distance can also be
investigated in the current setting. See the seminal work by~\mbox{\cite{azzalinibowmanhardle89,bickelrosenblatt93}}
and~\cite{hardlemammen93}. They are not qualitatively different as
shown in the classical goodness of fit tests.

Since the testing problem (\ref{2.17}) is equivalent to the following testing
problem:
\begin{equation} \label{2.20.1}
\cases{
H_{0}\dvtx P(z|x,2\Delta)-R(z|x,2\Delta)=0 &\quad\mbox{for all} $(x,z)\in
S^{2}$,
\cr
H_{1}\dvtx P(z|x,2\Delta)-R(z|x,2\Delta)\neq0 &\quad\mbox{for
some} $(x,z)\in S^{2}$,
}
\end{equation}
with, in light of (\ref{2.19}),
\[
R(z|x,2\Delta)=\int_{-\infty}^{z}r(t|x,2\Delta) \,dt=E\{P(z|Y,\Delta)|X=x\},
\]
then transition distribution-based tests can be formulated too. Let $\hat
{P}(z|x,2\Delta)$ be the direct estimator for the $2\Delta$-transition
distribution
\begin{equation}\label{3.7'}
\hat{P}(z|x,2\Delta)=\frac{1}{n}\sum_{i=1}^{n}W_{n}(X_{i}-x,x;h_{1}%
)I(Z_{i}<z).
\end{equation}
Regressing the transition distribution $P(z|X_{j},\Delta)$ on $X_{j-1}$ yields
$\hat{R}(z|x,2\Delta)$:
\begin{equation}\label{3.7''}
\hat{R}(z|x,2\Delta)=n^{-1}\sum_{i=1}^{n}W_{n}(X_{i}-x,x;h_{3})\hat{P}%
(z|Y_{i},\Delta),
\end{equation}
where $\hat{P}(z|y,\Delta)=n^{-1}\sum_{i=1}^{n}W_{n}(Y_{i}-y,y;b_{1}%
)I(Z_{i}<z).$ Similarly to (\ref{2.20}), for the testing problem
(\ref{2.20.1}), the transition distribution-based test will be
\begin{equation} \label{2.21}
T_{2}=\sum_{i=1}^{n}\{\hat{P}(Z_{i}|X_{i},2\Delta)-\hat{R}(Z_{i}|X_{i}%
,2\Delta)\}^{2}\omega(X_{i}),
\end{equation}
where the weight function $\omega(\cdot)$ is chosen to depend on only
$x$-variable, because $\hat{P}(z|x,2\Delta)$ is a nonparametric estimator of
the conditional distribution function, and we need only to weight down the
contribution from the sparse regions in the $x$-coordinate.

Note that the test statistic $T_{2}$ involves only one-dimensional smoothing.
Hence, it is expected to be more stable than $T_{1}$, and the null
distribution of $T_{2}$ can be better approximated by the asymptotic null
distribution. This will be justified by the theorems in the next section.

The choice between the transition density and distribution-based tests
reflects different degrees of smoothness of alternatives that we wish to test.
In a simpler problem of the traditional goodness-of-fit tests, this has been
thoroughly studied in~\cite{fan96}. Essentially, the transition density-based tests are more powerful in detecting local deviations whereas the
transition distribution-based tests are more powerful for detecting global deviations.

\section{Asymptotic properties}\label{sec:asy}

\subsection{Assumptions}

We assume the following conditions. These conditions are frequently imposed
for nonparametric studies for dependent data.

\renewcommand{\theassumption}{(A\arabic{assumption})}
\setcounter{assumption}{0}
\begin{assumption}
\label{ass:A1}The observed time series $\{X_{i}\}_{i=1}^{n+2}$ is strictly
stationary with time-homogenous $j\Delta$-transition density $p(X_{i+j}%
|X_{i},j\Delta).$
\end{assumption}

\begin{assumption}
\label{ass:A2}The kernel functions $W$ and $K$ are symmetric and bounded
densities with bounded supports, and satisfy the Lipschitz condition.
\end{assumption}

\begin{assumption}
\label{ass:A3}The weight function $w(x,z)$ has a continuous second-order
derivative with a compact support $\Omega^{\ast}$.
\end{assumption}

\begin{assumption}
\label{ass:A4}The stationary process $\{X_{i}\}$ is $\beta$-mixing with the
exponential decay rate $\beta(n)=O(e^{-\lambda n})$ for some $\lambda>0.$
\end{assumption}

\begin{assumption}
\label{ass:A5}The functions $p(y|x;\Delta)$ and $p(z|x;2\Delta)$ have
continuous second-order partial derivatives with respect to $(x,y)$ and
$(x,z)$ on the set $\Omega^{\ast}.$ The invariant density $\pi(x)$ of
$\{X_{i}\}$ has a continuous second-order derivative for $x\in\Omega_{x}%
^{\ast}$, a project of the set $\Omega^{\ast}$ onto the $x$-axis. Moreover,
$\pi(x)>0$, $p(y|x,\Delta)>0$ and $p(z|x,2\Delta)>0$ for all $(x,y)\in
\Omega^{\ast}$ and $(x,z)\in\Omega^{\ast}.$

\end{assumption}

\begin{assumption}
\label{ass:A6}The joint density $p_{1\ell}(x_{1},x_{\ell})$ of $(X_{1}%
,X_{\ell})$ for $\ell>1$ is bounded by a constant independent of $\ell$. Put
$g_{1\ell}(x_{1},x_{\ell})=p_{1\ell}(x_{1},x_{\ell})-\pi(x_{1})\pi(x_{\ell})$.
The function $g_{1\ell}$ satisfies the Lipschitz condition: for all
$(x^{\prime},y^{\prime})$ and $(x,y)$ in $\Omega^{\ast}$,
\[
|g_{1\ell}(x,y)-g_{1\ell}(x^{\prime},y^{\prime})|\leq C\sqrt{(x-x^{\prime
})^{2}+(y-y^{\prime})^{2}}.
\]
\end{assumption}

\begin{assumption}
\label{ass:A7}The bandwidths $h_{i}$s and $b_{i}$ are of the same order and
satisfy $nh_{1}^{3}/\log n\rightarrow\infty$ and $nh_{1}^{5}\rightarrow0.$

\end{assumption}

\begin{assumption}\label{ass:A8}
The bandwidth $h_{1}$ converges to zero in such a way that $nh_{1}%
^{9/2}\rightarrow0$ and $nh_{1}^{3/2}\rightarrow\infty.$

\end{assumption}

\subsection{Asymptotic null distributions}

To introduce our asymptotic results, we need the following notation. For any
integrable function $f(x)$, let $\|f\|^{2}=\int f^{2}(x)\, dx$
and
\[
s(z|x,2\Delta)=\int p^{2}(z|y,\Delta)p(y|x,\Delta)\, dy=E[  p^{2}%
(z|Y_{1},\Delta)|X_{1}=x]  .
\]
Note that the sampled observations $\{X_{n+2-i}\}_{i=0}^{n+1}$ are a reverse
Markov process under the null model. We also use $p^{\ast}(x|z,2\Delta)$ to
denote the $2\Delta$-transition density of the reverse process,
and let
\[
s^{\ast}(x|z,2\Delta)=\int p^{\ast2}(y|z,\Delta)p^{\ast}(x|y,\Delta)\, dy.
\]

Denote by
\begin{eqnarray*}
\Omega_{11}&=&\int w(x,z)p^{2}(z|x,2\Delta)\, dx\,dz,
\\
   \Omega_{12}&=&\int
w(x,z)p^{3}(z|x,2\Delta)\, dx\,dz,
\\
\Omega_{13}&=&\int w(x,z)s(z|x,2\Delta)p(z|x,2\Delta) \,dx\, dz,
\\
\Omega_{14}&=&\int w(x,z)r^{2}(z|x,2\Delta)p(z|x,2\Delta) \,dx\, dz,
\\
\Omega_{15}&=&\int w(x,z)s^{*}(x|z,2\Delta)p^{*}(x|z,2\Delta)[\pi(z)/\pi
(x)]^{2}\, dx\, dz,
\\
\Omega_{2}&=&\int w^{2}(x,z)p^{4}(z|x,2\Delta)\, dx\,dz.
\end{eqnarray*}
For a kernel function $K(\cdot)$, let $K^{*}(\cdot)=K*K(\cdot)$ and
$K_{h}(\cdot)=h^{-1}K(\cdot/h).$ Denote by $V(x,z)$ the conditional variance
function of $P(z|Y,\Delta)$, given $X=x$. Then it is easy to see that
\begin{eqnarray*}
\Omega_{13}-\Omega_{14}  &  =&\int w(x,z)V(x,z)p(z|x,2\Delta)\, dx\,dz
\\
&=&E\{V(X,Z)w(X,Z)|X=x\}.
\end{eqnarray*}
Throughout the paper, we use the notation $T_{n} \stackrel{a}{\sim} \chi
_{a_{n}}^{2}$ for a diverging sequence of constants $a_{n}$ to represent that
\[
(T_{n} - a_{n})/\sqrt{2 a_{n}} \stackrel{\mathcal{D}}{\longrightarrow
}{\mathcal{N}}(0,1).
\]

\begin{theorem}
\label{Th1} Assume Conditions \textup{\ref{ass:A1}--\ref{ass:A7}} hold. If $\{X_{i}\}$ is Markovian,
\[
(T_{1}-\mu_{1})/\sigma_{1}\stackrel{\mathcal{D}}{\longrightarrow}{\mathcal{N}%
}(0,1),
\]
where
\begin{eqnarray*}
\mu_{1}  &=&\Omega_{11}\|W\|^{2}\|K\|^{2}/(h_{1}h_{2})-\Omega_{12}%
\|W\|^{2}h_{1}^{-1}
\\
&&{}  +(\Omega_{13}-\Omega_{14})\|W\|^{2}/h_{3}+\Omega_{15}\|K\|^{2}/b_{2},
\end{eqnarray*}
and $\sigma_{1}^{2}=2\Omega_{2}\|W*W\|^{2}\|K*K\|^{2}/(h_{1}h_{2})$.
Furthermore, $r_{1}T_{1}\stackrel{a}{\sim}\chi^{2}_{a_{n}}$, where $a_{n}%
=r_{1}\mu_{1}$ and $r_{1}=2\mu_{1}/\sigma_{1}^{2}.$
\end{theorem}

The test statistic $T_{1}^{*}$, as far as its null distribution is concerned,
can be regarded as a special case of $T_{1}$, with the weight function
$w(x,z)=p^{-2}(z|x,2\Delta)w^{*}(x,z)$. Correspondingly, let $\Omega_{1j}^{*}$
denote $\Omega_{1j}$ with $w(x,z)$ replaced by $p^{-2}(z|x,2\Delta)\times w^{*}(x,z)$
and $\Omega_{2}^{*}$ defined similarly. Then, we have

\begin{corollary}
Under the conditions in Theorem \textup{\ref{Th1}} with $w$ replaced by $w^{*}$,
$r_{1}^{*}T_{1}^{*}\stackrel{a}{\sim}\chi^{2}_{a_{n}^{*}},$ where
\begin{eqnarray*}
r_{1}^{*}&=&\frac{\Omega_{11}^{*}\|W\|^{2}\|K\|^{2}}{\Omega_{2}^{*}%
\|W*W\|^{2}\|K*K\|^{2}}\bigl(1+o(1)\bigr),
\\
a_{n}^{*}&=&\frac{\Omega_{11}^{*}{}^{2}\|W\|^{4}\|K\|^{4}}{\Omega_{2}%
^{*}\|W*W\|^{2}\|K*K\|^{2}}\frac{1}{h_{1}h_{2}}\bigl(1+o(1)\bigr).
\end{eqnarray*}
\end{corollary}

The $r_{1}^{\ast}$ is asymptotically a constant depending on only the kernels
and the weight function. The degree of freedom $a_{n}^{\ast}$ is independent
of nuisance parameters. This reflects that the Wilks phenomenon continues to
hold in the current situation.

\begin{theorem}
\label{Th2} Under Conditions \textup{\ref{ass:A1}--\ref{ass:A6}} and \textup{\ref{ass:A8}}, if $\{X_{i}\}$ is Markovian,
\[
(T_{2}-\mu_{2})/\sigma_{2}\stackrel{\mathcal{D}}{\longrightarrow}{\mathcal{N}%
}(0,1),
\]
where
\[
\mu_{2}=\frac{1}{6h_{1}}\|W\|^{2}\int\omega(x) \{1+6h_{1}h_{3}^{-1}%
E[V(X_{\Delta},Z_{\Delta})|X_{\Delta}=x]\} \,dx,
\]
and $\sigma_{2}^{2}=\|W*W\|^{2} \|\omega\|^{2}/(45h_{1})$. Furthermore,
$r_{2}T_{2}\stackrel{a}{\sim}\chi^{2}_{b_{n}}$, where $b_{n}=r_{2}\mu_{2}$ and
$r_{2}=2\mu_{2}/\sigma_{2}^{2}$.
\end{theorem}

Comparing Theorems \ref{Th1} and \ref{Th2}, it is seen that asymptotic
variance of $T_{1}$ is an order of magnitude larger than that of $T_{2}$.
Therefore, the null distribution of $T_{2}$ can be more stably approximated
than that of $T_{1}$. On the other hand, the degrees of freedom in $T_{1}$ are
larger than in $T_{2}$, and the transition density-based tests are more
omnibus, capable of testing a wider class of alternative hypothesis.

\subsection{Power under contiguous alternative models}

To assess the power of the tests, we consider the following contiguous
alternative sequence for $T_{1}$:
\begin{eqnarray} \label{3b.1}
H_{1n}\dvtx  p(z|x,2\Delta)-r(z|x,2\Delta)=g_{n}(x,z),
\end{eqnarray}
where $g_{n}$ satisfies $E[g_{n}^{2}(X,Z)]=O(\delta_{n}^{2})$ and
$\operatorname{var}[g_{n}^{2}(X,Z)]\le M(E[g_{n}^{2}(X,Z)])^{2}$ for a constant
$M>0$ and a sequence $\delta_{n}$ going to zero as $n\to\infty$. Then the
power of the test statistic $T_{1}$ can be approximated using the following theorem.

\begin{theorem}
\label{Th3} Under Conditions \textup{\ref{ass:A1}--\ref{ass:A7}}, if $nh_{1}h_{2}\delta_{n}^{2}=O(1)$,
then under the alternative hypothesis $H_{1n}$,
\begin{eqnarray*}
(T_{1}-\mu_{1}-d_{1n})/\sigma_{1n}\stackrel{{\mathcal{D}}}{\longrightarrow
}{\mathcal{N}}(0,1),
\end{eqnarray*}
where $d_{1n}=nE\{g_{n}^{2}(X,Z)w(X,Z)\}(1+o(1))$, and $\sigma_{1n}%
=\sqrt{\sigma_{1}^{2}+4\sigma_{1A}^{2}}$ with
\[
\sigma_{1A}^{2}=nE[g_{n}^{2}(X,Z)w^{2}(X,Z)\{p(Z |X,2\Delta
)-p^{2}(Z|X,2\Delta)\}^{2}].
\]
\end{theorem}

Using Theorem \ref{Th1}, one can construct an approximate level-$\alpha$ test
based on $T_{1}$. Let $c_{\alpha}$ be the critical value such that
\[
P\{(T_{1}-\mu_{1})/\sigma_{1}\geq c_{\alpha}|H_{0}\}\leq\alpha.
\]
Then we have the following result, which demonstrates that the test statistic
$T_{1}$ can detect alternatives at rate $\delta_{n}=O(n^{-2/5})$.

\begin{theorem}
\label{Th4} Under Conditions \textup{\ref{ass:A1}--\ref{ass:A6}}, $T_{1}$ can detect alternatives with
rate $\delta_{n}=O(n^{-2/5})$ when $h_{1}=c_{1}n^{-1/5}$ and $h_{2}%
=c_{2}n^{-1/5}$ for some constants $c_{1}$ and $c_{2}$. Specifically, if
$\delta_{n}=dn^{-2/5}$ for a constant $d$, then:
\begin{eqnarray*}
&&\phantom{i}\mbox{\textup{(i)}}\quad\limsup_{d\to0}\limsup_{n\to\infty} P\{(T_{1}-\mu
_{1})/\sigma_{1}\ge c_{\alpha}|H_{1n}\}\le\alpha;
\\
&&\mbox{\textup{(ii)}}\quad\liminf_{d\to\infty}\liminf_{n\to\infty} P\{(T_{1}%
-\mu_{1})/\sigma_{1}\ge c_{\alpha}|H_{1n}\}=1.
\end{eqnarray*}
\end{theorem}

Similarly to (\ref{3b.1}), we consider the following alternative sequence to
study of the power of the test statistic $T_{2}$:
\[
H_{2n}\dvtx P(z|x,2\Delta)-R(z|x,2\Delta)=G_{n}(x,z),
\]
where $G_{n}(x,z)$ satisfies $E[G_{n}^{2}(X,Z)]=O(\rho_{n}^{2})$ and
$\operatorname{var}(G_{n}^{2}(X,Z))\le M\times \break (E[G_{n}^{2}(X,Z)])^{2}$ for a constant
$M>0$ and a sequence $\rho_{n}$ tending to zero. Then using the following
theorem one can calculate the power of the test statistic $T_{2}$.

\begin{theorem}
\label{Th5} Under Conditions \textup{\ref{ass:A1}--\ref{ass:A6}} and \textup{\ref{ass:A8}}, if $nh_{1}h_{3}\rho_{n}%
^{2}=O(1)$, then under the alternative hypothesis $H_{2n}$,
\begin{eqnarray*}
(T_{2}-\mu_{2}-d_{2n})/\sigma_{2n} \stackrel{\mathcal{D}}{\longrightarrow
}{\mathcal{N}}(0,1),
\end{eqnarray*}
where $d_{2n}=nE[G_{n}^{2}(X,Z)\omega(X)]+O(nh_{1}^{2}\rho_{n}+\rho_{n}%
h_{1}^{-1})$, $\sigma_{2n}^{2}=\sigma_{2}^{2}+4\sigma_{2A}^{2}$ and
\begin{eqnarray*}
\sigma_{2A}^{2}  &=&nE\biggl[\int G_{n}(X,Z)\omega(X )I(Z<z)P(dz|X,2\Delta)\biggr]^{2}
\\
&&{}-nE\biggl[\int G_{n}(X,Z)\omega(X)P(z|X,2\Delta)P(dz|X,2\Delta)\biggr]^{2}.
\end{eqnarray*}
\end{theorem}

In a manner parallel to Theorem \ref{Th4}, the following theorem demonstrates the
optimality of the test.

\begin{theorem}
\label{Th6} Under Conditions \textup{\ref{ass:A1}--\ref{ass:A6}}, $T_{2}$ can detect alternatives with
rate $\rho_{n}=O(n^{-4/9})$ when $h_{1}=c_{*}n^{-2/9}$ for some constant
$c_{*}$.
\end{theorem}

From Theorem \ref{Th6}, $T_{2}$ can detect alternatives at rate $O(n^{-4/9})$.
Using an argument similar to \cite{fanjiang05}, we can also establish the
minimax rate, $O(n^{-4/9})$, of the test. Note that the rate is optimal
according to \cite{ingster93,lepskispokoiny99} and \cite{spokoiny96}.
Compared with Theorem \ref{Th4}, it is seen that $T_{2}$ is more powerful than
$T_{1}$ for testing the Markov hypothesis. This is due to the fact that the
alternative under consideration for $T_{2}$ is global, namely, the density
under the alternative is basically globally shifted away from the null
hypothesis. On the other hand, $T_{1}$ and $T_{1}^{\ast}$ are more powerful
than $T_{2}$ for detecting local features of the alternative hypothesis. We
will now explore these features in simulations.

\section{Simulations}\label{sec:sim}

An important application of our test methods is to verify the Markov property
in the context where the null model is a diffusion process, since it is often
assumed in modern financial theory and practice that the observation process
comes from an underlying diffusion. Hence, we consider simulations for the
diffusion models.

To use the test statistics, one needs to find their null distributions.
Theoretically the asymptotic null distributions may be used to determine the
$p$-values of the test statistics. However, in practical applications the
asymptotic distributions do not necessarily give accurate approximations,
since the local sample size $nh_{1}h_{2}$ may not be large enough. This
phenomenon is shared by virtually all nonparametric kinds of tests where some
form of functional estimation is used.

We will mainly focus on the finite sample performance of the test statistic
$T_{1}^{\ast}$, since it possesses the Wilks property which facilitates
bandwidth selection and determination of the null distribution using a
bootstrap method.
Since the asymptotic null distribution of $T_{1}^{\ast}$ is independent of
nuisance parameters/functions under the null hypothesis, for a finite sample
it does not sensitively depend on the nuisance parameters/functions.
Therefore, the null distribution can be approximated by bootstraps, by fixing
nuisance parameters/functions at their reasonable estimates, as in
\cite{fanjiang07} in a different context.

In general, different bootstrap approximations to the null distributions are
needed for different null models, partially due to the large family of null
models with the Markov property. We will illustrate this method for the
Ornstein--Uhlenbeck model, which in financial mathematics is used for instance
as the \cite{vasicek77} model for interest rates. For other parametric models,
our approach can similarly be applied.

The Ornstein--Uhlenbeck model employed as the null hypothesis is
\begin{equation}\label{eq:OU}
dX_{t}=\kappa(\alpha-X_{t})\,dt+\sigma \,dW_{t},
\end{equation}
where $W_{t}$ is a Brownian motion, and the parameters are set as $\kappa=0.2,$
$\alpha=0.085,$ $\sigma=0.08$, which are realistic for interest rates over
long periods. We simulated the model $1000$ times. In each simulation, we
draw a sample with sample size $n=2400$ and weekly sampling interval
$\Delta=1/52$ using for this purpose a higher frequency Euler approximation,
or an exact discretization. The bandwidth selection for the test statistic
$T_{1}^{\ast}$ is performed using the simple empirical rule proposed by
\cite{hyndmanyao02}. Alternative methods include the cross-validation
approaches of \cite{fanyim04} and~\cite{hallracineli04}, but their computation
is intensive especially when repeated many times in Monte Carlo.

Given a sample from the model, we fit the model using the least squares method
and obtain the residuals of the fit, and then generate bootstrap samples using
the residual-based bootstrap method. For each simulation, we obtained three
bootstrap samples (this is merely for the reduction of computation cost; using
more samples will not fundamentally alter the results) and computed the test
statistic $T_{1}^{\ast}$ using the same bandwidths as the original sample
in the simulation. Pooling together the bootstrap samples from each
simulation, we obtained $3000$ bootstrap statistics. Their sampling
distributions, computed via the kernel density estimate, is taken as the
distribution of the bootstrap method. By using the kernel density estimation
method, the distribution of the realized values of the test statistic
$T_{1}^{\ast}$ in simulations is obtained as the true distribution (except for
the Monte Carlo errors).


\begin{figure}[b]

\includegraphics{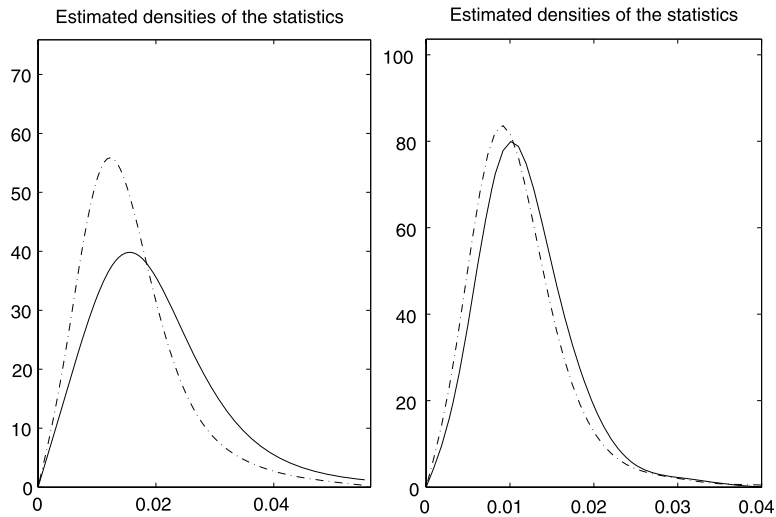}

\caption{Estimated densities. Left panel: $n=1200$; right panel: $n=2400$.
Solid---true, dotted---the bootstrap approximation.}
\label{fig1}
\end{figure}

Figure \ref{fig1} displays the estimated densities for $T_{1}^{\ast}$. Not
surprisingly, the bootstrapped distributions get much closer to the true ones
as the sample sizes increase. In our experience, the bootstrap approximations
start to become adequate for sample sizes starting at about $2400$.


\begin{table}[b]
\caption{Power of the test against $H_{1 \theta}$}
\label{tab1}
\begin{tabular*}{262pt}{@{\extracolsep{4in minus 4in}}lccccccc@{}}
\hline
&  & \multicolumn{6}{c@{}}{\textbf{Parameter} $\bolds{\theta}$}\\[-6pt]
&  & \multicolumn{6}{c@{}}{\hrulefill}\\
$\bolds{s}$ & \textbf{Level} $\bolds{\alpha}$ & \textbf{0.0} & \textbf{0.2} & \textbf{0.4} & \textbf{0.6} & \textbf{0.8} & \textbf{1.0}\\
\hline
$100$ & $0.01$ & $0.011$ & $1$\phantom{.000} & $1$\phantom{.000} & $1$\phantom{.000} & $1$\phantom{.000} & $1$\phantom{.000}\\
& $0.05$ & $0.055$ & $1$\phantom{.000} & $1$\phantom{.000} & $1$\phantom{.000} & $1$\phantom{.000} & $1$\phantom{.000}\\[3pt]
\phantom{0}$10$ & $0.01$ & $0.011$ & $0.010$ & $0.070$ & $0.228$ & $0.580$ & $0.846$\\
& $0.05$ & $0.055$ & $0.019$ & $0.123$ & $0.549$ & $0.901$ & $0.989$\\
\hline
\end{tabular*}
\end{table}


%


To investigate the power of the test statistics, we employ various sequences
of alternatives indexed by a parameter $\theta=0,$ $0.2,$ $0.4,$ $0.6,$ $0.8,$
$1.0$. One of the main ways for an otherwise Markovian model to become
non-Markovian is to restrict too much its state space. For instance, consider
a bivariate diffusion model. Taken jointly, the two components are Markovian,
but taken in isolation a single component may not be:

\begin{enumerate}[1.]
\item Alternative model with missing state variable in the drift: we first
consider the situation where the null model (\ref{eq:OU}) is missing a state
variable, in this case $X$ mean-revers to the stochastic level $\theta
\alpha_{t}+(1-\theta)\alpha$ under the alternative
\[
H_{1\theta}\dvtx dX_{t}=\kappa\bigl(  \theta\alpha_{t}+(1-\theta)\alpha
-X_{t}\bigr) \,  dt+\sigma \,dW_{t},
\]
where $\alpha_{t}$ is the random process
\[
d\alpha_{t}=\kappa_{1}(a-\alpha_{t}) \,dt+\sigma_{1}\, dB_{t},
\]
with $B_{t}$ a the Brownian motion independent of $W_{t}$, $\kappa_{1}%
=\kappa/s$, $a=s\alpha$, and $\sigma_{1}=\sigma/2$, with $s=100$ and $10$.
When $\theta\neq0$, the alternatives are non-Markovian. The results in the
first part of Table~\ref{tab1} show that the test statistic rejects the null
hypothesis when the observations are drawn under $H_{1\theta}.$

\begin{table}[b]
\caption{Power of the test against $H_{2 \theta}$}
\label{tab2}
\begin{tabular*}{262pt}{@{\extracolsep{4in minus 4in}}lccccccc@{}}
\hline
&  & \multicolumn{6}{c@{}}{\textbf{Parameter} $\bolds{\theta}$}\\[-6pt]
&  & \multicolumn{6}{c@{}}{\hrulefill}\\
$\bolds{s}$ & \textbf{Level} $\bolds{\alpha}$ & \textbf{0.0} & \textbf{0.2} & \textbf{0.4} & \textbf{0.6} & \textbf{0.8} & \textbf{1.0}\\
\hline
$1000$ & $0.01$ & $0.013$ & $0.402$ & $0.660$ & $0.762$ & $0.813$ & $0.817$\\
& $0.05$ & $0.067$ & $0.557$ & $0.768$ & $0.845$ & $0.878$ & $0.905$\\[3pt]
\phantom{0}$100$ & $0.01$ & $0.013$ & $0.028$ & $0.183$ & $0.372$ & $0.492$ & $0.573$\\
& $0.05$ & $0.067$ & $0.098$ & $0.340$ & $0.527$ & $0.627$ & $0.697$\\[3pt]
\phantom{00}$10$ & $0.01$ & $0.013$ & $0.007$ & $0.020$ & $0.017$ & $0.032$ & $0.088$\\
& $0.05$ & $0.067$ & $0.037$ & $0.052$ & $0.070$ & $0.122$ & $0.218$\\
\hline
\end{tabular*}
\end{table}

\item Alternative model with missing state variable in volatility: next, we
consider alternative models where volatility is stochastic,
\[
H_{2\theta}\dvtx dX_{t}=\kappa(\alpha-X_{t})\, dt+\bigl((1-\theta)\sigma+\theta\sigma
_{t}\bigr)\, dW_{t},
\]
where $\sigma_{t}=\sqrt{Y_{t}}$ is a random process following the \cite{cir85}
model
\[
dY_{t}=\kappa_{2}(b-Y_{t})\, dt+\sigma_{2}Y_{t}^{1/2}\, dB_{2t},
\]
where $B_{2t}$ is a standard Brownian motion independent of $W_{t}$,
$\kappa_{2}=\kappa/s$, $b=s\alpha$ and $\sigma_{2}=\sigma/2$, with $s=1000$,
$100$ and $10$. When $\theta\neq0$, the alternatives are also non-Markovian.

\item Alternative model with missing state variable in jumps: finally, we
consider a model with compound Poisson jumps
\[
H_{3\theta}\dvtx dX_{t}=\kappa(\alpha-X_{t})\, dt+\sigma \,dW_{t}+J_{t}
 \,dN_{t}(\theta),
\]
where $N_{t}(\theta)$ is a Poisson process with stochastic intensity $\theta$
and jump size $1$, while $J_{t}$ is a the jump size. We will consider two
types of jump sizes:

\begin{enumerate}[(ii)]
\item[(i)] $J_{t}$ is independent of ${\mathcal{F}}_{t}$ and follows
$N(0,\sigma_{1}^{2})$ with $\sigma_{1}=\sigma/2$, which makes $H_{3\theta}$ Markovian;

\item[(ii)] $J_{t}$ follows the CIR model
\[
dJ_{t}=\kappa(a-J_{t})\, dt+\sigma_{1}J_{t}^{1/2}\,dB_{3t},
\]
where $B_{3t}$ is a standard Brownian motion independent of $W_{t}$, $K=0.2$,
$a=0.085$ and $\sigma_{1}=0.08/2$. Then $J_{t}$ is not independent of
${\mathcal{F}}_{t}$. This leads to alternatives $H_{3\theta}$ which are not
Markovian for $\theta\neq0$.
\end{enumerate}
\end{enumerate}

The alternative models considered here are $\beta$-mixing.  For example, in the first alternative
$H_{1\theta}$, the joint process $(X_t, \alpha_t)$ is an affine process and it is $\beta$-mixing.
Hence, $X_t$ is $\beta$-mixing.  A similar argument applies to two other alternatives.
In fact, for the first alternative $H_{1\theta}$, the time series $(X_{i\Delta}, \alpha_{i\Delta})$
can be written as a bivariate autoregressive model.  Hence, it is $\beta$-mixing with the choice of parameters.
Note that for all of the above alternatives, when $\theta$ is small, the null
and alternative models are nearly impossible to differentiate. In the limit
where $\theta=0$, the null and the alternative are identical. Therefore, it
can be expected that, when $\theta=0$, the power of test should be close to
the significance level; and as $\theta$ deviates more from $0$, the power
should increase. Also we can expect that our tests will be able to detect only
the type (ii) jumps but not the type (i) jump, since for the type (i) jump the
alternatives are Markovian.

The simulated powers are reported in Tables~\ref{tab1}--\ref{tab3}.
The null distribution of the normalized test statistics does not depend sensitively on choice of
bandwidth, whereas the power depends on the choice of bandwidth and the alternative under consideration.
As expected, our test is fairly powerful for detecting non-Markovian alternatives
$H_{k\theta}$ $(k=1,2,3),$ at least in situations where the alternative is
sufficiently far from the null. For $H_{3\theta}$, the test has, as it should,
no power to identify the type (i) alternatives but is powerful for
discriminating against the type (ii) alternatives. This illustrates well the
sensitivity and specificity of our tests.

%

\begin{table}[t]
\caption{Power of the test against $H_{3 \theta}$}
\label{tab3}
\begin{tabular*}{287pt}{@{\extracolsep{4in minus 4in}}lccccccc@{}}
\hline
& &\multicolumn{6}{c@{}}{\textbf{Parameter} $\bolds{\theta}$}  \\[-6pt]
& &\multicolumn{6}{c@{}}{\hrulefill}\\
\textbf{Jump type} & \textbf{Level} $\bolds{\alpha}$ & \textbf{0.0} & \textbf{0.2} & \textbf{0.4} &\textbf{0.6} & \textbf{0.8} & \textbf{1.0}\\
\hline
\phantom{i}(i) & $0.01$ & $0.010$ & $0.009$ & $0.023$ & $0.003$ & $0.016$ & $0.009$\\
& $0.05$ & $0.059$ & $0.048$ & $0.054$ & $0.054$ & $0.058$ & $0.056$\\[3pt]
(ii) & $0.01$ & $0.010$ & $0.514$ & $0.774$ & $0.888$ & $0.940$ & $0.951$\\
& $0.05$ & $0.059$ & $0.533$ & $0.796$ & $0.894$ & $0.946$ & $0.961$\\
\hline
\end{tabular*}
\end{table}


\eject
\section{Technical proofs}\label{sec:tech}

\subsection{Technical lemmas}

We now introduce some technical lemmas, the proofs of which can be found in the
supplemental material of this paper. To save space, some notation in the
lemmas will appear later in the course of proofs of the main theorems.


\begin{lemma}
\label{le0} Suppose that $W$ is symmetric and continuous with a bounded
support. If $h\to0$ and $nh\to\infty$, then
\begin{eqnarray*}
W_{n}(z,x;h)  & =& \biggl\{\frac{1}{\mu_{0}(W)\pi(x)}-\frac{z}{h}\frac{h\pi^{\prime
}(x)}{\pi^{2}(x)\mu_{0}(W)} +O_{p}(\rho_{n}(h))\biggr\}W_{h}(z)
\\
&&{}  +O_{p}(\rho_{n}(h))\frac{z}{h}W_{h}(z),
\end{eqnarray*}
uniformly for $x\in\Omega^{*}$, where $O_{p}(\rho_{n}(h))$ does not depend on
$z$, where $\mu_{0}(W) = \int W(u)\, du$.
\end{lemma}

\begin{lemma}
\label{le1} Under Conditions \textup{\ref{ass:A1}--\ref{ass:A6}}:
\begin{longlist}[(iii)]
\item[(i)] for $k=0,1,$
\[
\sup_{(y,z)\in\Omega^{*}}\Biggl|n^{-1}\sum_{i=1}^{n}b_{1}^{-k}(Y_{i}-y)^{k}W_{b_{1}%
}(Y_{i}-y)\varepsilon_{i}(z)\Biggr| =O\bigl\{\sqrt{(\log n)/(nb_{1}b_{2})}\bigr\};
\]

\item[(ii)] for $k=0,1$,
\[
\sup_{(x,z)\in\Omega^{*}}\Biggl|\frac{1}{n}\sum_{j=1}^{n}h_{3}^{-k}%
(X_{j}-x)^{k} W_{h_{3}}(X_{j}-x)e_{j}(z)\Biggr|
=O_{p}\bigl\{\sqrt{\log(n)/(nh_{3})}\bigr\};
\]

\item[(iii)] $\sup_{(x,z)\in\Omega^{*}}|\frac{1}{n}\sum_{j=1}^{n}%
q^{*}(x,Z_{j})\varepsilon_{j+1}(z)|
=O_{p}\{\sqrt{\log(n)/(nb_{2})}\};$\vspace*{1pt}

\item[(iv)] $\sup_{(x,z)\in\Omega^{*}}|\frac{1}{n}\sum_{j=1}^{n}W_{h_{1}%
}(X_{j}-x)\varepsilon_{j}^{*}(z)| =O_{p}\{\sqrt{(\log n)/(nh_{1}h_{2})}\}$.
\end{longlist}
\end{lemma}

\begin{lemma}
\label{le2} Under Conditions \textup{\ref{ass:A1}--\ref{ass:A6}}, we have
\[
\xi_{n}(x,y)\equiv\frac{1}{n}\sum_{i=1}^{n} r_{n1}(x,Y_{i})\varepsilon_{i}(z)
=O_{p}\bigl(\sqrt{n^{-1}b_{1}\log n} \bigr),
\]
uniformly for $(x,z)\in\Omega^{*}$, where $r_{n1}$ is defined right after
\textup{(\ref{ap9b})}.
\end{lemma}

\begin{lemma}
\label{le3} Suppose  Conditions \textup{\ref{ass:A1}--\ref{ass:A5}} hold. Then
\[
\eta_{n}(x,z)\equiv n^{-1}\sum_{i=1}^{n} r_{n}^{*}(x,Y_{i})\varepsilon_{i}(z)
=O\bigl\{\sqrt{[(b_{1}^{4}+h_{3}^{4})\log n]/(nb_{2})} \bigr\},
\]
uniformly for $(x,y)\in\Omega^{*}$, where $r_{n}^{*}(\cdot,\cdot)$ is defined
in \textup{(\ref{a.15a})}.
\end{lemma}

\begin{lemma}
\label{lead1} Under Conditions \textup{\ref{ass:A1}--\ref{ass:A6}}:
\begin{longlist}[(ii)]
\item[(i)] $\sum_{1\le i<j\le n}[\tilde{\psi}(i,j)-\tilde{\psi}(i)-\tilde
{\psi}(j)+\tilde{\psi}(0)] =o_{p}(h_{1}^{-1})$;

\item[(ii)] $(n-1)\sum_{i=1}^{n}[\tilde{\psi}(i)-\tilde{\psi}(0)]=o_{p}%
(h_{1}^{-1}).$
\end{longlist}
\end{lemma}

\begin{lemma}
\label{lead2} Assume Conditions \textup{\ref{ass:A1}--\ref{ass:A5}} hold. Then we have:
\begin{longlist}[(ii)]
\item[(i)] under Condition \textup{\ref{ass:A6}},
\begin{eqnarray*}
\tfrac{1}{2}n(n-1)\tilde{\psi}(0)  & =&\Omega_{11}\|W\|^{2}\|K\|^{2}%
/(h_{1}h_{2})-\Omega_{12}\|W\|^{2}/h_{1}
\\
&&{}  +\Omega_{13}\|W\|^{2}/h_{3} -\Omega_{14}\|W\|^{2}/h_{3}
\\
&&{}  +\Omega_{15}\|K\|^{2}/b_{2}+O(n^{-2});
\end{eqnarray*}

\item[(ii)] under Condition \textup{\ref{ass:A7}},
\begin{eqnarray*}
&&\frac{1}{2}n(n-1)\tilde{\phi}(0)
\\
 &&\qquad =\frac{1}{6h_{1}}\|W\|^{2}\int\omega(x)
  \{1+6h_{1}h_{3}^{-1}E[V(X_{\Delta},Z_{\Delta})|X_{\Delta
}=x]\} \,dx
+O(1).
\end{eqnarray*}
\end{longlist}
\end{lemma}

\begin{lemma}
\label{lead3} Assume that Conditions \textup{\ref{ass:A1}--\ref{ass:A5}} hold. Then we have:

\begin{longlist}[(ii)]
\item[(i)] under Condition \textup{\ref{ass:A6}},
\[
\sigma_{1n}^{-1}\sum_{1\le i<j\le n}\psi^{*}(i,j)\stackrel{\mathcal{D}%
}{\longrightarrow}{\mathcal{N}}(0,1),
\]
where $\sigma_{1n}^{2}=2\Omega_{2}\|W*W\|^{2}\|K*K\|^{2}/(n^{2}h_{1}h_{2});$

\item[(ii)] under Condition \textup{\ref{ass:A7}},
\[
\sigma_{2n}^{-1}\sum_{1\le i<j\le n}\phi^{*}(i,j)\stackrel{\mathcal{D}%
}{\longrightarrow}{\mathcal{N}}(0,1),
\]
where $\sigma_{2n}^{2}=\|W*W\|^{2}\|w\|^{2}/(45n^{2}h_{1}).$
\end{longlist}
\end{lemma}

\subsection{Preliminaries}

Since the test statistics $T_{1}$ and $T_{1}^{*}$ compare the difference
between $\hat{p}(z|x,2\Delta)$ and $\hat{r}(z|x,2\Delta)$, we derive an
asymptotic expression for this difference under $H_{0}$ before giving the
proofs of theorems. In addition, in order to streamline our arguments, we will
introduce some technical lemmas and put them behind the proofs of theorems.
The arguments employed here use techniques from the $U$-statistic and
nonparametric smoothing.

First let us introduce some notation. Let $\rho_{n}(h)=h^{2}+\sqrt{\log
n/(nh)},$\break $\mu_{0}(W)=\int W(x)\, dx$ and $\mu_{2}(W)=\int x^{2}W(x) \,dx$.
Denote by $m(y,z)=\break E\{K_{b_{2}}(Z_{j}-z)|Y_{j}=y\}$, $m^{*}(x,z)=E\{K_{h_{2}%
}(Z_{j}-z)|X_{j}=x\}$, $m_{1}(y,z)=\break\partial m(y,z)/\partial y$ and $m_{1}%
^{*}(x,z)=\partial m^{*}(x,z)/\partial x$.

Using an elementary property of the local linear smoother (see, e.g.,
\cite{fanyao03}), we obtain that
\begin{equation} \label{ap1}
\hat{p}(z|x,2\Delta)-p(z|x,2\Delta)=A_{n}^{\ast}(x,z)+B_{n}^{\ast}%
(x,z)+C_{n}^{\ast}(x,z),
\end{equation}
where $\varepsilon_{j}^{\ast}(z)=K_{h_{2}}(Z_{j}-z)-m^{\ast}(X_{j},z),$
\begin{eqnarray}\label{ap1a0}
A_{n}^{\ast}(x,z)&=&\frac{1}{n}\sum_{j=1}^{n}W_{n}(X_{j}-x,x;h_{1}%
)\varepsilon_{j}^{\ast}(z),\nonumber
\\
B_{n}^{\ast}(x,z)&=&\frac{1}{n}\sum_{j=1}^{n}W_{n}(X_{j}-x,x;h_{1})\nonumber
\\[-8pt]\\[-8pt]
&&\hspace*{25pt}{}\times \{m^{\ast
}(X_{j},z)-m^{\ast}(x,z)-m_{1}^{\ast}(x,z)(X_{j}-x)\},\nonumber
\\
C_{n}^{\ast}(x,z)&=&m^{\ast}(x,z)-p(z|x,2\Delta).\nonumber
\end{eqnarray}
By
a second-order Taylor expansion,
\[
B_{n}^{\ast}(x,z)=\frac{1}{n}\sum_{j=1}^{n}W_{n}(X_{j}-x,x;h_{1})\frac
{h_{1}^{2}}{2}m_{2}^{\ast}(\tilde{x}_{j},z)\biggl(\frac{X_{j}-x}{h_{1}%
}\biggr)^{2},
\]
where $m_{2}^{\ast}(\tilde{x},z)=\frac{\partial^{2}m^{\ast}(x,z)}{\partial
^{2}x}|_{x=\tilde{x}_{j}},$ and $\tilde{x}_{j}$ lies between $X_{j}$ and
$x$.
By \cite{fanyaotong96}, it is easy to show that
\begin{equation} \label{ap1a}
B_{n}^{\ast}(x,z)=O_{p}(h_{1}^{2}) \quad \mbox{and} \quad C_{n}^{\ast
}(x,z)=O_{p}(h_{2}^{2}),
\end{equation}
uniformly for $(x,z)\in\Omega^{\ast}.$ By the definition of $\hat{r}$, we
have
\begin{equation}\label{ap2}
\hat{r}(z|x,2\Delta)-r(z|x,2\Delta)=L_{n1}(x,z)+L_{n1}^{\ast}(x,z),
\end{equation}
where
\begin{eqnarray*}
L_{n1}(x,z)&=&\frac{1}{n}\sum_{j=1}^{n}W_{n}(X_{j}-x,x;h_{3})\{\hat{p}%
(z|Y_{j},\Delta)-p(z|Y_{j},\Delta)\},
\\
L_{n1}^{\ast}(x,z)&=&\frac{1}{n}\sum_{j=1}^{n}W_{n}(X_{j}-x,x;h_{3}%
)\{p(z|Y_{j},\Delta)-r(z|x,2\Delta)\}.
\end{eqnarray*}
Subtracting (\ref{ap2}) from (\ref{ap1}), we obtain that, under $H_{0}%
\dvtx p(z|x,2\Delta)=r(z|x,2\Delta),$
\begin{eqnarray}\label{ap3}
 \quad \hat{p}(z|x,2\Delta)-\hat{r}(z|x,2\Delta)&=& A_{n}^{\ast}(x,z)+B_{n}^{\ast
}(x,z)  +C_{n}^{\ast}(x,z)\nonumber
\\[-8pt]\\[-8pt]
&&{}-L_{n1}(x,z) -L_{n2}(x,z)-L_{n3}(x,z),\nonumber
\end{eqnarray}
where
\begin{eqnarray*}
L_{n2}(x,z)&=&n^{-1}\sum_{j=1}^{n}W_{n}(X_{j}-x,x;h_{3})\{p(z|Y_{j}%
,\Delta)-r(z|X_{j},2\Delta)\},
\\
L_{n3}(x,z)&=&n^{-1}\sum_{j=1}^{n}W_{n}(X_{j}-x,x;h_{3})\{r(z|X_{j}%
,2\Delta)-r(z|x,2\Delta)\}.
\end{eqnarray*}
By the continuity of $\partial^{2}r(z|x,2\Delta)/\partial x^{2}$, it is easy
to show that
\begin{equation} \label{ap4}
L_{n3}(x,z)=O_{p}(h_{3}^{2}) \qquad \mbox{uniformly for }  (x,z)\in
\Omega^{\ast}.
\end{equation}
Therefore, by (\ref{ap1a}), (\ref{ap3}) and (\ref{ap4}),
\begin{eqnarray}\label{ap4a}
&&\hat{p}(z|x,2\Delta)-\hat{r}(z|x,2\Delta)\nonumber
\\[-8pt]\\[-8pt]
&&{}\qquad =[A_{n}^{\ast}(x,z)-L_{n2}(x,z)]
-L_{n1}(x,z)+O_{p}\Biggl(\sum_{i=1}^{3}h_{i}^{2}\Biggr).\nonumber
\end{eqnarray}
Let $e_{j}(z)=p(z|Y_{j},\Delta)-r(z|X_{j},2\Delta)$. Then it can be rewritten
that
\begin{equation}\label{ap5}
L_{n2}(x,z)=\frac{1}{n}\sum_{j=1}^{n}W_{n}(X_{j}-x,x;h_{3})e_{j}(z).
\end{equation}
Note that $r(z|X_{j},2\Delta)=E\{p(z|Y_{j},\Delta)|X_{j}\}$. It follows that
$E[e_{j}(z)|X_{j}]=0$ and $\operatorname{Var}[e_{j}(z)]=O(1)$ uniformly for
$z$ and $j=1,\ldots,n$. Applying Lemma \ref{le0} with $z=X_{j}-x$ and $h=h_{3}%
$, we obtain that
\begin{eqnarray}\label{ap9}
&&W_{n}(X_{j}-x,x;h_{3})\nonumber
\\
&&\qquad =\biggl\{\frac{1}{\mu_{0}\pi(x)}-\frac{X_{j}%
-x}{h_{3}}\frac{h_{3}\pi^{\prime}(x)}{\pi^{2}(x)\mu_{0}}+O_{p}(\rho_{n}%
(h_{3}))\biggr\}W_{h_{3}}(X_{j}-x)
\\
&&{}\qquad\quad +O_{p}(\rho_{n}(h_{3}))\frac{X_{j}-x}{h_{3}}W_{h_{3}}(X_{j}-x),\nonumber
\end{eqnarray}
uniformly for $x\in\Omega^{\ast}$, where $O_{p}(\rho_{n}(h_{3}))$ does not
depend on $j$. Therefore,
\[
L_{n2}(x,z)=L_{n21}(x,z)-L_{n22}(x,z)+L_{n23}(x,z)+L_{n24}(x,z),
\]

\noindent
where
\begin{eqnarray*}
L_{n21}(x,z)&=&\frac{1}{\mu_{0}(W)\pi(x)}n^{-1}\sum_{j=1}^{n}W_{h_{3}}%
(X_{j}-x)e_{j}(z),
\\
L_{n22}(x,z)&=&\frac{h_{3}\pi^{\prime}(x)}{\mu_{0}(W)\pi^{2}(x)}n^{-1}\sum
_{j=1}^{n}\frac{X_{j}-x}{h_{3}}W_{h_{3}}(X_{j}-x)e_{j}(z),
\\
L_{n23}(x,z)&=&O_{p}(\rho_{n}(h_{3}))n^{-1}\sum_{j=1}^{n}W_{h_{3}}(X_{j}%
-x)e_{j}(z),
\\
L_{n24}(x,z)&=&O_{p}(\rho_{n}(h_{3}))n^{-1}\sum_{j=1}^{n}\frac{X_{j}-x}{h_{3}%
}W_{h_{3}}(X_{j}-x)e_{j}(z).
\end{eqnarray*}
By Lemma \ref{le1}(ii), we have $L_{n21}(x,z)=O_{p}\{\sqrt{(\log n)/(nh_{3}%
)}\}$ and
\[
L_{n22}(x,z)=O_{p}\bigl\{h_{3}\sqrt{(\log n)/(nh_{3})}\bigr\}=O_{p}\bigl\{\sqrt{(h_{3}\log
n)/n}\bigr\},
\]
uniformly for $(x,z)\in\Omega^{\ast}$.
Then
\[
L_{n23}(x,z)=O_{p}\bigl\{\rho_{n}(h_{3})\sqrt{(\log n)/(nh_{3})}\bigr\}=o_{p}%
\bigl\{\sqrt{(h_{3}\log n)/n}\bigr\}
\]
and $L_{n24}(x,z)=o_{p}\{\sqrt{(h_{3}\log n)/n}\}$, uniformly for
$(x,z)\in\Omega^{\ast}$. Then
\[
L_{n2}(x,z)=\frac{1}{\mu_{0}(W)\pi(x)}\frac{1}{n}\sum_{j=1}^{n}W_{h_{3}}%
(X_{j}-x)e_{j}(z)+O_{p}\bigl\{\sqrt{(h_{3}\log n)/n}\bigr\},
\]
uniformly for $(x,z)\in\Omega^{\ast}$.
Note that from (\ref{ap1a0}) and (\ref{ap5})
\begin{eqnarray} \label{ap6a}
\quad A_{n}^{\ast}(x,z)-L_{n2}(x,z)  &=&\frac{1}{n}\sum_{j=1}^{n}[W_{n}%
(X_{j}-x,x;h_{1})\varepsilon_{j}^{\ast}(z)
\\
&&\qquad\hspace*{4pt}{}  -W_{n}(Y_{j}-x,x;h_{3})e_{j+1}(z)]+r_{n}(x,z),
\end{eqnarray}
where
\[
r_{n}(x,z)
=-\frac{1}{n}W_{n}(X_{1}-x,x;h_{3})e_{1}(z)+\frac{1}{n}W_{n}(Y_{n}%
-x,x;h_{3})e_{n+1}(z),
\]
which is of order $O_{p}(1/(nh_{3}))=o_{p}\{\sqrt{(h_{3}\log n)/{n}}\}$,
uniformly for $(x,z)\in\Omega^{\ast}$. Let $\varepsilon_{i}(z)=K_{b_{2}}%
(Z_{i}-z)-m(Y_{i},z).$ Then, similarly to (\ref{ap1}), we have
\begin{equation}\label{ap7}
\hat{p}(z|y,\Delta)-p(z|y,\Delta)=A_{n}(y,z)+B_{n}(y,z)+C_{n}(y,z),
\end{equation}
where $A_{n}(y,z)=n^{-1}\sum_{i=1}^{n}W_{n}(Y_{i}-y,y;b_{1})\varepsilon
_{i}(z)$, $B_{n}(y,z)=O_{p}(b_{1}^{2})$ and\break $C_{n}(y,z)=O_{p}(b_{1}^{2})$,
uniformly for $(y,z)\in\Omega^{\ast}$. It follows from the definition of
$L_{n1}$ that
\begin{eqnarray}\label{ap8}
L_{n1}(x,z)  &  =&n^{-1}\sum_{i=1}^{n}W_{n}(X_{j}-x,x;h_{3})A_{n}%
(Y_{j},z)\nonumber
\\[-8pt]\\[-8pt]
&&{}  +n^{-1}\sum_{i=1}^{n}W_{n}(X_{j}-x,x;h_{3})[B_{n}(Y_{j},z)+C_{n}%
(Y_{j},z)].\nonumber
\end{eqnarray}
Using Lemma \ref{le0}, we get
\[
A_{n}(y,z)=A_{n1}(y,z)-A_{n2}(y,z)+A_{n3}(y,z)+A_{n4}(y,z),
\]
where
\begin{eqnarray*}
A_{n1}(y,z)&=&\frac{1}{\mu_{0}\pi(y)}n^{-1}\sum_{i=1}^{n}W_{b_{1}}%
(Y_{i}-y)\varepsilon_{i}(z),
\\[-2pt]
A_{n2}(y,z)&=&\frac{b_{1}\pi^{\prime}(y)}{\mu_{0}\pi^{2}(x)}n^{-1}\sum_{i=1}%
^{n}\frac{Y_{i}-y}{b_{1}}W_{b_{1}}(Y_{i}-y)\varepsilon_{i}(z),
\\[-2pt]
A_{n3}(y,z)&=&O_{p}(\rho_{n}(b_{1}))n^{-1}\sum_{i=1}^{n}W_{b_{1}}(Y_{i}%
-y)\varepsilon_{i}(z),
\\[-2pt]
A_{n4}(y,z)&=&O_{p}(\rho_{n}(b_{1}))n^{-1}\sum_{i=1}^{n}\frac{Y_{i}-y}{b_{1}%
}W_{b_{1}}(Y_{i}-y)\varepsilon_{i}(z).
\end{eqnarray*}
Using Lemma \ref{le1}(i), we obtain that
\[
A_{n3}(y,z)=O_{p}(\rho_{n}(b_{1}))O_{p}\biggl(\sqrt{\frac{\log n}{nb_{1}b_{2}}%
}\biggr)\vspace*{-2pt}
\]
and\vspace*{-2pt}
\[
A_{n4}(y,z)=O_{p}(\rho_{n}(b_{1}))O_{p}\biggl(\sqrt{\frac{\log n}{nb_{1}b_{2}}%
}\biggr),
\]
uniformly for $(y,z)\in\Omega^{\ast}$. Then
\[
A_{n}(y,z)=A_{n1}(y,z)-A_{n2}(y,z)+O_{p}(\rho_{n}(b_{1}))\sqrt{(\log
n)/(nb_{1}b_{2})},
\]
uniformly for $(y,z)\in\Omega^{\ast}$. This, combined with (\ref{ap8}) and
Condition \ref{ass:A6}, yields that
\begin{equation} \label{ap9a}
\qquad L_{n1}(x,z) =L_{n11}(x,z)-L_{n12}(x,z)+L_{n13}(x,z) +O_{p}\{(\log
n)/(nb_{1}^{3/2})\},
\end{equation}
where $L_{n11}(x,z)=n^{-1}\sum_{j=1}^{n}W_{n}(X_{j}-x,x;h_{3})A_{n1}%
(Y_{j},z)$,
\begin{eqnarray*}
L_{n12}(x,z)&=&n^{-1}\sum_{j=1}^{n}W_{n}(X_{j}-x,x;h_{3})A_{n2}(Y_{j},z),
\\[-2pt]
L_{n13}(x,z)&=&n^{-1}\sum_{j=1}^{n}W_{n}(X_{j}-x,x;h_{3})[B_{n}(Y_{j}%
,z)+C_{n}(Y_{j},z)].
\end{eqnarray*}

\noindent
Note that, by Lemma \ref{le1}(i), $A_{n1}(y,z)=O_{p}\{\sqrt{(\log
n)/(nb_{1}b_{2})} \}$, uniformly for $(y,z)\in\Omega^{\ast}.$ Using Lemma
\ref{le0}, we obtain that
\begin{equation}\label{ap9b}
L_{n11}(x,z)=M_{n11}(x,z)+M_{n12}(x,z)+O_{p}\{(\log n)/(nb_{1}^{3/2})\},
\end{equation}
where
\begin{eqnarray*}
M_{n11}(x,z)&=&\frac{1}{\mu_{0}^{2}\pi(x)}\frac{1}{n^{2}}\sum_{j=1}^{n}%
\sum_{i=1}^{n}W_{h_{3}}(X_{j}-x)W_{b_{1}}(Y_{i}-Y_{j})\pi^{-1}(Y_{j}%
)\varepsilon_{i}(z),
\\
M_{n12}(x,z)&=&\frac{h_{3}\pi^{\prime}(x)}{\mu_{0}^{2}\pi^{2}(x)}\frac{1}{n^{2}%
}\sum_{j=1}^{n}\sum_{i=1}^{n}\frac{X_{j}-x}{h_{3}}W_{h_{3}}(X_{j}-x)W_{b_{1}%
}(Y_{i}-Y_{j})
\\
&&{}\hspace*{86pt}\times \pi^{-1}(Y_{j})\varepsilon_{i}(z).
\end{eqnarray*}
Let
\[
M_{n11}^{\ast}(x,y)=n^{-1}\sum_{j=1}^{n}W_{h_{3}}(X_{j}-x)W_{b_{1}}%
(y-Y_{j})\pi^{-1}(Y_{j}),
\]
$g_{n}(x,y)=E[M_{n11}^{\ast}(x,y)]$ and $r_{n1}(x,y)=M_{n11}^{\ast
}(x,y)-g_{n}(x,y)$. Then
\[
M_{n11}(x,z)=\frac{1}{\mu_{0}^{2}\pi(x)}\frac{1}{n}\sum_{i=1}^{n}g_{n}%
(x,Y_{i})\varepsilon_{i}(z)+\frac{1}{\mu_{0}^{2}\pi(x)}\frac{1}{n}\sum
_{i=1}^{n}r_{n1}(x,Y_{i})\varepsilon_{i}(z).
\]
By Lemma \ref{le2},
\begin{equation}\label{ap11}
M_{n11}(x,z)=\frac{1}{\mu_{0}^{2}\pi(x)}\frac{1}{n}\sum_{i=1}^{n}g_{n}%
(x,Y_{i})\varepsilon_{i}(z)+O_{p}\bigl\{\sqrt{(b_{1}\log n)/n} \bigr\},
\end{equation}
uniformly for $(x,z)\in\Omega^{\ast}.$ Similarly to Lemma \ref{le1}(iii), the
first term on the right-hand side of (\ref{ap11}) is $O_{p}\{\sqrt{(\log
n)/(nb_{2})}\}$, uniformly for $(x,z)\in\Omega^{\ast}$. Hence,
\[
\sup_{(x,z)\in\Omega^{\ast}}|M_{n11}(x,z)|=O_{p}\bigl\{\sqrt{(\log n)/(nb_{2}%
)} \bigr\}.
\]
Similarly, we have
\[
\sup_{(x,z)\in\Omega^{\ast}}|M_{n12}(x,z)|=O_{p}\bigl\{h_{3}\sqrt{(\log
n)/(nb_{2})} \bigr\}=O_{p}\bigl\{\sqrt{(b_{1}\log n)/n} \bigr\}.
\]
By the symmetry of the kernel function and Taylor's expansion, it can be shown
that
\begin{eqnarray*}
g_{n}(x,y)  &=&E[\pi^{-1}(Y_{1})W_{b_{1}}(y-Y_{1})W_{h_{3}}%
(X_{1}-x)]
\\
&=&\mu_{0}^{2}p(y|x,\Delta)\pi(x)/\pi(y)+O(b_{1}^{2}+h_{3}^{2})
\\
&\equiv& \mu_{0}^{2}p^{\ast}(x|y,\Delta)+O(b_{1}^{2}+h_{3}^{2}),
\end{eqnarray*}
uniformly for $(x,y)\in\Omega^{\ast}$, where $p^{\ast}(x|y,\Delta)$ is the
one-$\Delta$ transition density of the reverse series $\{X_{n+2-i}%
\}_{i=1}^{n+1}$, that is, the conditional density of $X_{1}$ given $Y_{1}=y$. Note
that $g_{n}$ is a deterministic function. It follows that
\begin{equation}\label{a.15a}
g_{n}(x,Y_{i})=\mu_{0}^{2}p^{\ast}(x|Y_{i},\Delta)+r_{n}^{\ast}(x,Y_{i}),
\end{equation}
where $r_{n}^{\ast}(x,Y_{i})$ is $\sigma(Y_{i})$-measurable and is of order
$O(b_{1}^{2}+h_{3}^{2})$ for $(x,Y_{i})\in\Omega^{\ast}$. This combined with
(\ref{ap11}) leads to
\begin{eqnarray}\label{ap11add}
L_{n11}(x,z)  &=& \frac{1}{n}\sum_{i=1}^{n}q^{\ast}(x,Y_{i})\varepsilon
_{i}(z)+\frac{O(1)}{n}\sum_{i=1}^{n}r_{n}^{\ast}(x,Y_{i})\varepsilon_{i}(z)\nonumber
\\[-8pt]\\[-8pt]
&&{}  +O_{p}\bigl(\{\log n/(nb_{1}^{3/2})\}+\{b_{1}(\log n)/n\}^{1/2}\bigr),\nonumber
\end{eqnarray}
where $q^{\ast}(x,y)=p(y|x,\Delta)/\pi(y)$. The first term in (\ref{ap11add})
is obviously
\[
\frac{1}{n}\sum_{i=1}^{n}q^{\ast}(x,Z_{i})\varepsilon_{i+1}(z)+O_{p}\biggl(\frac
{1}{nb_{1}}\biggr).
\]
By Lemma \ref{le3}, the second term in (\ref{ap11add}) is $O_{p}(\sqrt
{(b_{1}^{4}+h_{3}^{4})\log(n)/(nb_{2})})$, uniformly for $(x,z)\in\Omega
^{\ast}.$ Then uniformly for $(x,z)\in\Omega^{\ast}$,
\[
L_{n11}(x,z)=\frac{1}{n}\sum_{i=1}^{n}q^{\ast}(x,Z_{i})\varepsilon
_{i+1}(z)+O_{p}\bigl(\{\log n/(nb_{1}^{3/2})\}+\{b_{1}(\log n)/n\}^{1/2}\bigr).
\]
In the same argument, $L_{n12}(x,z)$ is dominated by $L_{n11}(x,z)$ and is of
order
\[
b_{1}L_{n11}(x,z)=O_{p}\bigl(\{\log n/(nb_{1}^{3/2})\}+\{b_{1}\log n/n\}^{1/2}\bigr),
\]
which combined with (\ref{ap9a}) leads to
\begin{eqnarray}\label{A.17b}
L_{n1}(x,z)  &=& \frac{1}{n}\sum_{i=1}^{n}q^{\ast}(x,Z_{i})\varepsilon
_{i+1}(z)+L_{n13}(x,z)\nonumber
\\[-8pt]\\[-8pt]
&&{}  +O_{p}\bigl(\{\log n/(nb_{1}^{3/2})\}+\{b_{1}\log n/n\}^{1/2}\bigr),\nonumber
\end{eqnarray}
uniformly for $(x,z)\in\Omega^{\ast}$. This together with (\ref{ap4a}) and
(\ref{ap6a}) yields the following asymptotic expression:
\begin{equation}\label{ap10}
\qquad \hat{p}(z|x,2\Delta)-\hat{r}(z|x,2\Delta) =T_{n1}(x,z)+T_{n2}(x,z)
+T_{n3}(x,z)+T_{n4}(x,z),
\end{equation}
where
\begin{eqnarray*}
T_{n1}(x,z)  &=& \frac{1}{n}\sum_{j=1}^{n}[W_{n}(X_{j}-x,x;h_{1}%
)\varepsilon_{j}^{\ast}(z)
\\
&&\hspace*{26pt}{}  -W_{n}(Y_{j}-x,x;h_{3})e_{j+1}(z)-q^{\ast}(x,Z_{j})\varepsilon
_{j+1}(z)],
\\
T_{n2}(x,z)&=&n^{-1}\sum_{j=1}^{n}W_{n}(X_{j}-x,x;h_{3})[B_{n}(Y_{j}%
,z)+C_{n}(Y_{j},z)],
\\
T_{n3}(x,z)&=&B_{n}^{\ast}(x,z)+C_{n}^{\ast}(x,z)+L_{n3}(x,z),
\\
T_{n4}(x,z)&=&O_{p}\bigl(\{\log n/(nh_{1}^{3/2})\}+\{b_{1}\log n/n\}^{1/2}+\{\log
n/(nb_{1}^{3/2})\}\bigr),
\end{eqnarray*}
uniformly for $(x,z)\in\Omega^{\ast}$.

\subsection{Proofs of theorems}

We now give the proofs of our main results.

\begin{pf*}{Proof of Theorem \ref{Th1}} (i) \textit{Approximate $T_{1}$ by a
$U$-statistic.} Let $w_{i}=w(X_{i},Z_{i})$. By (\ref{ap10}) and the definition
of $T_{1}$, we have
\begin{eqnarray*}
T_{1}  & =&\sum_{i=1}^{n}w_{i}[T_{n1}(X_{i},Z_{i})+T_{n2}(X_{i},Z_{i}%
)+T_{n3}(X_{i},Z_{i})+T_{n4}(X_{i},Z_{i})]^{2}
\\
& =&\sum_{i=1}^{n}\sum_{k=1}^{4}w_{i}T_{nk}^{2}(X_{i},Z_{i})+2\sum_{i=1}%
^{n}w_{i}T_{n1}(X_{i},Z_{i})T_{n2}(X_{i},Z_{i})
\\
&&{}  +2\sum_{i=1}^{n}w_{i}T_{n1}(X_{i},Z_{i})T_{n3}(X_{i},Z_{i})+2\sum_{i=1}%
^{n}w_{i}T_{n2}(X_{i},Y_{i})T_{n3}(X_{i},Y_{i})
\\
&&{}  +2\sum_{i=1}^{n}w_{i}[T_{n1}(X_{i},Z_{i})+T_{n2}(X_{i},Z_{i})+T_{n3}%
(X_{i},Z_{i})]T_{n4}(X_{i},Z_{i})
\\
& \equiv& T_{11}+T_{12}+T_{13}+T_{14}+T_{15}.
\end{eqnarray*}
By Lemmas \ref{le0} and \ref{le1}, $T_{n1}(x,z)=O_{p}\{\sqrt{(\log
n)/(nh_{1}h_{2})}\}$. Note that $T_{n2}(x,z)=O_{p}(b_{1}^{2})$, $T_{n3}%
(x,z)=O_{p}(h_{1}^{2})$, uniformly for $(x,z)\in\Omega^{\ast}$. It is
straightforward to verify that $T_{14}=O_{p}(nh_{1}^{4})=o(1/h_{1})$,
$T_{15}=o_{p}(1/\sqrt{h_{1}h_{2}})$. Using the same argument as for (B.2) in \cite{yacfanpeng09},
we obtain $T_{12}=o_{p}(1/\sqrt{h_{1}h_{2}})$ and $T_{13}=o_{p}(1/\sqrt
{h_{1}h_{2}})$. Therefore,
\[
T_{1}=\sum_{i=1}^{n}\sum_{k=1}^{4}w_{i}T_{nk}^{2}(X_{i},Z_{i})+o_{p}%
(h_{1}^{-1}).
\]
Note that
\begin{eqnarray*}
\sum_{i=1}^{n}w_{i}T_{n2}^{2}(X_{i},Z_{i})&=&O_{p}(nh_{1}^{4})=o_{p}(1/h_{1}),
\\
\sum_{i=1}^{n}w_{i}T_{n3}^{2}(X_{i},Z_{i})&=&o_{p}(1/h_{1})
\end{eqnarray*}
and
\[
\sum_{i=1}^{n}w_{i}T_{n4}^{2}(X_{i},Z_{i})=o_{p}(1/h_{1})
.
\]
It follows that
\begin{eqnarray*}
T_{1}  &  =&\sum_{i=1}^{n}w_{i}T_{n1}^{2}(X_{i},Z_{i})+o_{p}(h_{1}^{-1})
\\
&  \equiv&\tilde{T}_{1}+o_{p}(h_{1}^{-1}).
\end{eqnarray*}
It can be rewritten that
\[
\tilde{T}_{1}=\sum_{i=1}^{n}w_{i}[B_{n1}^{\ast}(X_{i},Z_{i})-B_{n2}^{\ast
}(X_{i},Z_{i})-B_{n3}(X_{i},Z_{i})]^{2},
\]
where $B_{n1}^{\ast}(x,z)=\frac{1}{n}\sum_{j=1}^{n}W_{n}(X_{j}-x,x;h_{1}%
)\varepsilon_{j}^{\ast}(z),$
\[
B_{n2}^{\ast}(x,z)=\frac{1}{n}\sum_{j=1}^{n}W_{n}(Y_{j}-x,x;h_{3})e_{j+1}(z)
\]
and
\begin{eqnarray*}
B_{n3}(x,z)  &  =&\frac{1}{n}\sum_{j=1}^{n}q^{\ast}(x,Z_{j})\varepsilon
_{j+1}(z)
\\
&  =&\frac{1}{n}\frac{1}{\pi(x)}\sum_{j=1}^{n}p(Z_{j}|x,\Delta)\pi(x)\pi
^{-1}(Z_{j})\varepsilon_{j+1}(z).
\end{eqnarray*}
Applying Lemmas \ref{le0} and \ref{le1} and using Condition \ref{ass:A5}, we obtain
that
\[
\tilde{T}_{1}=\sum_{i=1}^{n}w_{i}\{B_{n1}(X_{i},Z_{i})-B_{n2}(X_{i}%
,Z_{i})-B_{n3}(X_{i},Z_{i})\}^{2}+o_{p}(h_{1}^{-1}),
\]
where $B_{n1}(x,z)=\frac{1}{n}\frac{1}{\pi(x)}\sum_{j=1}^{n}W_{h_{1}}%
(X_{j}-x)\varepsilon_{j}^{\ast}(z)$ and
\[
B_{n2}(x,z)=\frac{1}{n}\frac{1}{\pi(x)}\sum_{j=1}^{n}W_{h_{3}}(Y_{j}%
-x)e_{j+1}(z).
\]
Hence,
\[
T_{1}=\sum_{i=1}^{n}w_{i}\{B_{n1}(X_{i},Z_{i})-B_{n2}(X_{i},Z_{i}%
)-B_{n3}(X_{i},Z_{i})\}^{2}+o_{p}(h_{1}^{-1}).
\]
Let
$\xi(i,j)=W_{h_{1}}(X_{j}-X_{i})\varepsilon_{j}^{\ast}(Z_{i})-W_{h_{3}}%
(Y_{j}-X_{i})e_{j+1}(Z_{i})-q(X_{i},Z_{j})\varepsilon_{j+1}(Z_{i})$ and
\[
\psi(i,j,k)={n^{-2}}{w_{i}}\pi^{-2}(X_{i})\xi(i,j)\xi(i,k),
\]
where $q(x,z)=p(z|x,\Delta)\pi(x)/\pi(z)=p^{\ast}(x|z,\Delta).$ Then
\[
T_{1}=\sum_{i,j,k=1}^{n}\psi(i,j,k)+o_{p}(h_{1}^{-1}).
\]

(ii) \textit{Derive the asymptotics using the asymptotic theory for the
U-statistic.} Let
\begin{eqnarray*}
B_{11}&=&\sum_{i<j<k}\{\psi(i,j,k)+\psi(i,k,j)+\psi(j,i,k)
\\
&&\hspace*{26pt}{}+\psi(j,k,i)+\psi
(k,i,j)+\psi(k,j,i)\},
\\
B_{12}&=&\sum_{i\neq j}[\psi(i,j,j)+\psi(j,i,j)+\psi(j,j,i)]
\end{eqnarray*}
and
\[
 B_{13}=\sum_{i=1}^{n}\psi(i,i,i).
 \]
Then
\begin{equation}\label{A.19}
T_{1}=B_{11}+B_{12}+B_{13}+o_{p}(h_{1}^{-1}).
\end{equation}
Let $\psi^{\ast}(i,j,k)=\psi(i,j,k)+\psi(i,k,j)+\psi(j,i,k)+\psi
(j,k,i)+\psi(k,i,j)+\psi(k,j,i)$. Then $\psi^{\ast}(i,j,k)$ is symmetrical
about $(i,j,k)$, and hence $B_{11}=\sum_{i<j<k}\psi^{\ast}(i,j,k).$ Using
Hoeffding's decomposition, we obtain that
\begin{equation}\label{newad1}
B_{11}=\sum_{i<j<k}\Phi(i,j,k)+(n-2)\sum_{1\leq i<j\leq n}\psi^{\ast}(i,j),
\end{equation}
where
\[
\Phi(i,j,k)=\psi^{\ast}(i,j,k)-\psi^{\ast}(i,j)-\psi^{\ast}(i,k)-\psi^{\ast
}(j,k),
\]
$\psi^{\ast}(i,j)=\int\psi^{\ast}(i,j,k) \,dF(x_{k},y_{k},z_{k})$ and $F$ is
the distribution of $(X_{k},Y_{k},Z_{k})$. Applying the lemma with
$\delta=1/3$ in \cite{gaoking04},
we can show that $E\{\sum_{i<j<k}\Phi(i,j,\break k)\}^{2}=o(h_{1}^{-2}).$ Therefore,
the first term on the right-hand side of (\ref{newad1}) is $o_{p}(h_{1}^{-1}%
)$, so that
\begin{equation}\label{newad1a}
B_{11}=(n-2)\sum_{1\leq i<j\leq n}\psi^{\ast}(i,j)+o_{p}(h_{1}^{-1}).
\end{equation}
By the Markovian property of $\{X_{i}\}$, $E[\psi^{\ast}(i,j)]=0$. Hence, up
to a\vspace*{1pt} ignorable term of order $o_{p}(h_{1}^{-1})$, $B_{11}$ is a
$U$-statistic with mean zero. Define $\tilde{\psi}(i,j)=\psi(i,i,j)+\psi
(i,j,i)+\psi(j,i,i)+\psi(j,j,i)+\psi(j,i,j)+\psi(i,j,j)$,\vspace*{1pt} $\tilde{\psi
}(i)=\int\tilde{\psi}(i,j)\, dF(x_{j},y_{j},z_{j})$ and $\tilde{\psi
}(0)=E[\tilde{\psi}(i)]$. Then we have
\[
B_{12}=\sum_{1\leq i<j\leq n}\tilde{\psi}(i,j).
\]
Since $\tilde{\psi}(i,j)$ is a symmetrical kernel, using the Hoeffding
decomposition, we obtain that
\begin{eqnarray}
B_{12}  &=&\sum_{1\leq i<j\leq n}[\tilde{\psi}(i,j)-\tilde{\psi}%
(i)-\tilde{\psi}(j)+\tilde{\psi}(0)]\nonumber
\\[-8pt]\\[-8pt]
&&{}  +(n-1)\sum_{i=1}^{n}[\tilde{\psi}(i)-\tilde{\psi}(0)]+\frac{1}%
{2}n(n-1)\tilde{\psi}(0).\nonumber
\end{eqnarray}
By Lemma \ref{lead1},
\begin{equation}
B_{12}=\tfrac{1}{2}n(n-1)\tilde{\psi}(0)+o_{p}(h_{1}^{-1}).
\end{equation}
Note that $B_{13}\geq0$. By straightforward calculation on the mean of
$B_{13}$, it can be shown that
\begin{equation}\label{newad3}
B_{13}=O_{p}\bigl(n/(n^{2}h_{1}^{2}h_{2}^{2})\bigr)=o_{p}(h_{1}^{-1}).
\end{equation}
Therefore, a combination of (\ref{A.19}) and (\ref{newad1a})--(\ref{newad3})
leads to
\begin{equation}\label{newad4}
T_{1}=\frac{1}{2}n(n-1)\tilde{\psi}(0)+(n-2)\sum_{1\leq i<j\leq n}\psi^{\ast
}(i,j)+o_{p}(h_{1}^{-1}).
\end{equation}
By Lemma \ref{lead2}(i),
\[
\tfrac{1}{2}n(n-1)\tilde{\psi}(0)=\mu_{1}+o_{p}(h_{1}^{-1}).
\]
Applying Lemma \ref{lead3}(i), we obtain that
\[
(n-2)\sum_{i<j}\psi^{\ast}(i,j)/\sigma_{1}\stackrel{\mathcal{D}}%
{\longrightarrow}{\mathcal{N}}(0,1),
\]
where $\sigma_{1}^{2}=2\Omega_{2}\|W\ast W\|^{2}\|K\ast K\|^{2}/(h_{1}h_{2}).$
Therefore, the result of this theorem holds.
\end{pf*}

\begin{pf*}{Proof of Theorem \ref{Th2}} The proof is similar to that of Theorem~\ref{Th1}.\vspace*{1pt}

(i) \textit{Asymptotic expression for $\hat{P}(z|x,2\Delta)-\hat
{R}(z|x,2\Delta)$}. By the definitions in (\ref{3.7'}) and (\ref{3.7''}),
\begin{eqnarray} \label{a1}
\hat{P}(z|x,2\Delta)-P(z|x,2\Delta)  &=& \frac{1}{n}\sum_{i=1}^{n} W_{n}%
(X_{i}-x,x;h_{1})\nonumber
\\[-8pt]\\[-8pt]
&&{}\qquad\ \times \lbrack I(Z_{i}<z)-P(z|x,2\Delta)],\nonumber
\\\label{a2}
\hat{R}(z|x,2\Delta)-R(z|x,2\Delta)&=&S_{n1}(x,z)+S_{n2}(x,z),
\end{eqnarray}
where $S_{n1}(x,z)=n^{-1}\sum_{i=1}^{n}W_{n}(X_{i}-x,x;h_{3})[\hat{P}%
(z|Y_{i},\Delta)-P(z|Y_{i},\Delta)]$ and
\begin{equation}
S_{n2}(x,z)=n^{-1}\sum_{i=1}^{n}W_{n}(X_{i}-x,x;h_{3})[P(z|Y_{i}%
,\Delta)-R(z|x,2\Delta)].
\end{equation}
Let $u_{i}(z,\Delta)=I(Z_{i}<z)-P(z|Y_{i},\Delta)$. Then $E[u_{i}%
(z,\Delta)]=0$. By (\ref{2.12}),
\[
\hat{P}(z|y,\Delta)-P(z|y,\Delta)=n^{-1}\sum_{i=1}^{n}W_{n}(Y_{i}%
-y,y;b_{1})[I(Z_{i}<z)-P(z|y,\Delta)].
\]
This can be rewritten as
\begin{equation}
\hat{P}(z|y,\Delta)-P(z|y,\Delta)=P_{n1}(y,z)+P_{n2}(y,z),
\end{equation}
where
\begin{eqnarray*}
P_{n1}(y,z)  &  =&n^{-1}\sum_{i=1}^{n}W_{n}(Y_{i}-y,y;b_{1})u_{i}(z,\Delta),
\\
P_{n2}(y,z)  &  =&n^{-1}\sum_{i=1}^{n}W_{n}(Y_{i}-y,y;b_{1})[P(z|Y_{i}%
,\Delta)-P(z|y,\Delta)].
\end{eqnarray*}
By Lemma \ref{le0} and the symmetry of the kernel function $W(\cdot)$, and by
using Taylor's expansion, it is easy to show that
\begin{equation}\label{Th2:aa2}
P_{n2}(y,z)=(\partial^{2}/\partial y^{2})P(z|y,\Delta)b_{1}^{2}+o_{p}%
(b_{1}^{2})=O_{p}(b_{1}^{2}),
\end{equation}
uniformly for $(y,z)\in\Omega^{\ast}.$ Hence,
\begin{equation}\label{A.275}
\hat{P}(z|y,\Delta)-P(z|y,\Delta)=P_{n1}(y,z)+O_{p}(b_{1}^{2}),
\end{equation}
uniformly for $(y,z)\in\Omega^{\ast}.$ Then
\begin{equation}
S_{n1}(x,z)=n^{-1}\sum_{i=1}^{n}W_{n}(X_{i}-x,x;h_{3})P_{n1}(Y_{i}%
,z)+O_{p}(b_{1}^{2}),
\end{equation}
uniformly for $(x,z)\in\Omega^{\ast}.$ Using the same arguments as those for
$L_{n11}(x,z)$ between (\ref{ap9b}) and (\ref{ap10}), we obtain that
\begin{eqnarray}\label{Th2:eq2}
S_{n1}(x,z)  & =&\frac{1}{n}\sum_{i=1}^{n}q^{\ast}(x,Y_{i})u_{i}%
(z,\Delta)\nonumber
\\
&&{}  +O_{p}\bigl(\{\log n/(nb_{1}^{3/2})\}+\{b_{1}(\log n)/n\}^{1/2}\bigr)\nonumber
\\[-8pt]\\[-8pt]
&=&\frac{1}{n}\sum_{i=1}^{n}q^{\ast}(x,Z_{i})u_{i+1}(z,\Delta)\nonumber
\\
&&{}  +O_{p}\biggl(  \frac{\log n}{nb_{1}^{3/2}}+\{b_{1}(\log n)/n\}^{1/2}\biggr).\nonumber
\end{eqnarray}
Rewrite $S_{n2}(x,z)$ as
\[
S_{n2}(x,z)=S_{n21}(x,z)+S_{n22}(x,z),
\]
where
\begin{eqnarray*}
S_{n21}(x,z)&=&n^{-1}\sum_{i=1}^{n}W_{n}(X_{i}-x,x;h_{3})[P(z|Y_{i}%
,\Delta)-R(z|X_{i},2\Delta)],
\\
S_{n22}(x,z)&=&n^{-1}\sum_{i=1}^{n}W_{n}(X_{i}-x,x;h_{3})[R(z|X_{i}%
,2\Delta)-R(z|x,2\Delta)].
\end{eqnarray*}
By the continuity of $\partial^{2}R(z|x,2\Delta)/\partial x^{2}$ and the same
argument as that for (\ref{Th2:aa2}), $S_{n22}(x,z)=O_{p}(h_{3}^{2}),$
uniformly for $(x,z)\in\Omega^{\ast}$. Let $e_{i}^{\ast}(z)=P(z|Y_{i}%
,\Delta)-R(z|X_{i},2\Delta)$. Then $E[e_{i}^{\ast}(z)|X_{i}]=0$, and
\begin{eqnarray}\label{a3}
S_{n2}(x,z)  &=&n^{-1}\sum_{i=1}^{n}W_{n}(X_{i}-x,x;h_{3})e_{i}^{\ast
}(z)+O_{p}(h_{3}^{2})\nonumber
\\[-8pt]\\[-8pt]
&=&n^{-1}\sum_{i=1}^{n}W_{n}(Y_{i}-x,x;h_{3})e_{i+1}^{\ast}(z)+O_{p}%
(h_{3}^{2}).\nonumber
\end{eqnarray}
By (\ref{a1}) and (\ref{a2}), under $H_{0}$, we have
\begin{equation} \label{Th2:eq1}
\hat{P}(z|x,2\Delta)-\hat{R}(z|x,2\Delta)=-S_{n1}(x,z)-S_{n2}(x,z)+S_{n3}%
(x,z),
\end{equation}
where, with $u_{j}^{\ast}(z,2\Delta)=I(Z_{j}<z)-P(z|X_{j},2\Delta)$,
\begin{eqnarray*}
S_{n3}(x,z)  & =& \frac{1}{n}\sum_{i=1}^{n}W_{n}(X_{i}-x,x;h_{1})[I(Z_{i}%
<z)-P(z|x,2\Delta)]
\\
&=&\frac{1}{n}\sum_{i=1}^{n}W_{n}(X_{i}-x,x;h_{1})u_{i}^{\ast}(z,2\Delta)
\\
&&{}  +\frac{1}{n}\sum_{i=1}^{n}W_{n}(X_{i}-x,x;h_{1})[P(z|X_{i},2\Delta
)-P(z|x,2\Delta)].
\end{eqnarray*}
Similarly to (\ref{Th2:aa2}), the second term above is of order $O_{p}(h_{1}%
^{2})$,
\begin{equation} \label{Th2:eq1a}
S_{n3}(x,z)=\frac{1}{n}\sum_{i=1}^{n}W_{n}(X_{i}-x,x;h_{1})u_{i}^{\ast
}(z,2\Delta)+O_{p}(h_{1}^{2}),
\end{equation}
uniformly for $(x,z)\in\Omega^{\ast}$. A combination of (\ref{Th2:eq2}%
)--(\ref{Th2:eq1a}) yields that
\begin{equation}\label{Th2:eq3}
\hat{P}(z|x,2\Delta)-\hat{R}(z|x,2\Delta)=T_{n1}^{\ast}(x,z)+T_{n2}^{\ast
}(x,z)+T_{n3}^{\ast}(x,z),
\end{equation}
where
\begin{eqnarray*}
T_{n1}^{\ast}(x,z)  & =&\frac{1}{n}\sum_{j=1}^{n}[W_{n}(X_{j}-x,x;h_{1}%
)u_{j}^{\ast}(z,2\Delta)
\\
&&{}\hspace*{27pt}  -W_{n}(Y_{j}-x,x;h_{3})e_{j+1}^{\ast}(z)-q^{\ast}(x,Z_{j})u_{j+1}%
(z,\Delta)],
\\
T_{n2}^{\ast}(x,z)&=&O_{p}(b_{1}^{2}+h_{1}^{2}+h_{3}^{2}%
),
\end{eqnarray*}
uniformly for  $(x,z)\in\Omega^{\ast}$,   and
\begin{eqnarray*}
T_{n3}^{\ast}(x,z)=O_{p}\bigl(\{\log n/(nb_{1}^{3/2})\}+\{b_{1}(\log n)/n\}^{1/2}%
\bigr),
\end{eqnarray*}
uniformly for  $(x,z)\in\Omega^{\ast}$.

(ii) \textit{Asymptotic normality of $T_{2}$}. Similar to (\ref{newad4}), we
have
\begin{equation}\label{Th2:eq4}
T_{2}=\frac{1}{2}n(n-1)\tilde{\phi}(0)+(n-2)\sum_{1\le i<j\le n}\phi
^{*}(i,j)+o_{p}(h^{-1}),
\end{equation}
where $\tilde{\phi}(0)$ and $\phi^{*}(i,j)$ are defined the same as
$\tilde{\psi}(0)$ and $\psi^{*}(i,j)$, respectively, but with $\psi$ replaced
by
\begin{eqnarray*}
\phi(i,j,k)   ={n^{2}}{w_{i}}\pi^{-2}(X_{i}) \eta(i,j)\eta(i,k),
\end{eqnarray*}
where
\begin{eqnarray*}
\eta(i,j)&=&W_{h_{1}}(X_{j}-X_{i})u_{j}^{*}(Z_{i},2\Delta)-W_{h_{3}}(Y_{j}%
-X_{i})e^{*}_{j+1}(Z_{i})
\\
&&{}-q(X_{i},Z_{j})u_{j+1}(Z_{i},\Delta).
\end{eqnarray*}
By Lemma \ref{lead2}(ii), we have
\begin{eqnarray}\label{Th2:eq5}
\tfrac{1}{2}n(n-1)\tilde{\phi}(0)=\mu_{2}+o_{p}(h_{1}^{-1}).
\end{eqnarray}
By Lemma \ref{lead3}(ii), we have
\begin{eqnarray}\label{Th2:eq6}
(n-2)\sum_{i<j}\phi^{*}(i,j)/\sigma_{2}\stackrel{\mathcal{D}}{\longrightarrow
}{\mathcal{N}}(0,1).
\end{eqnarray}
A combination of (\ref{Th2:eq4})--(\ref{Th2:eq6}) completes the proof of the theorem.
\end{pf*}

\begin{pf*}{Proof of Theorem \ref{Th3}} Under $H_{1n}$, $p(z|x,2\Delta
)=r(z|x,2\Delta)+g_{n}(x,z).$ Similarly to (\ref{ap3}), we have under $H_{1n}$
\[
\hat{p}(z|x,2\Delta)-\hat{r}(z|x,2\Delta)=Q_{n}(x,z)+g_{n}(x,z),
\]
where
\[
Q_{n}(x,z)=A_{n}^{\ast}(x,z)+B_{n}^{\ast}(x,z)+C_{n}^{\ast}(x,z)-L_{n1}%
(x,z)-L_{n2}(x,z)-L_{n3}(x,z).
\]
Then
\begin{eqnarray}\label{Th3:eq1}
T_{1}  &=&\sum_{i=1}^{n}Q_{n}^{2}(X_{i},Z_{i})w_{i}+\sum_{i=1}^{n}g_{n}%
^{2}(X_{i},Z_{i})w_{i}\nonumber
\\[-8pt]\\[-8pt]
&&{}  +2\sum_{i=1}^{n}g_{n}(X_{i},Z_{i})Q_{n}(X_{i},Z_{i})w_{i}.\nonumber
\end{eqnarray}
Since $\delta_{n}^{2}=O(\frac{1}{nh_{1}h_{2}})$, it can be shown that
\begin{equation}\label{Th3:eq2}
\sum_{i=1}^{n}g_{n}^{2}(X_{i},Z_{i})w_{i}=nE[g_{n}^{2}%
(X,Z)w(X,Z)]+o_{p}\bigl(1/\sqrt{h_{1}h_{2}}\bigr).
\end{equation}
By (\ref{ap1a}) and (\ref{ap4}), $B_{n}^{\ast}(x,z)=O_{p}(h_{1}^{2})$,
$C_{n}^{\ast}(x,z)=O_{p}(h_{2}^{2})$ and $L_{n3}(x,z)=O_{p}(h_{3}^{2})$,
uniformly for $(x,z)\in\Omega^{\ast}$. It follows from the H\"{o}lder inequality
that
\begin{eqnarray}\label{a.34b}
&&  2\sum_{i=1}^{n}w_{i}g_{n}(X_{i},Z_{i})[B_{n}^{\ast}(X_{i},Z_{i}%
)+C_{n}^{\ast}(X_{i},Z_{i})-L_{n3}(X_{i},Z_{i})]\nonumber
\\[-8pt]\\[-8pt]
&&\qquad   =O_{p}\bigl(n\delta_{n}(h_{1}^{2}+h_{2}^{2}+h_{3}^{2})\bigr).\nonumber
\end{eqnarray}
A combination of (\ref{Th3:eq1})--(\ref{a.34b}) yields that
\begin{eqnarray}\label{Th3:eq3}
T_{1}  &=& \sum_{i=1}^{n}Q_{n}^{2}(X_{i},Z_{i})w_{i}+nE[g_{n}%
^{2}(X,Z)w(X,Z)]\nonumber
\\
&&{}  +2\sum_{i=1}^{n}g_{n}(X_{i},Z_{i})w_{i}[A_{n}^{\ast}(X_{i},Z_{i}%
)-L_{n2}(X_{i},Z_{i})-L_{n1}(X_{i},Z_{i})]\nonumber
\\[-8pt]\\[-8pt]
&&{}  +\bigl\{o_{p}\bigl(1/\sqrt{h_{1}h_{2}}\bigr)+O_{p}\bigl(n\delta_{n}(h_{1}^{2}+h_{2}^{2}%
+h_{3}^{2})\bigr)\bigr\}\nonumber
\\
&\equiv& T_{11}+T_{12}+T_{13}+o_{p}\bigl(1/\sqrt{h_{1}h_{2}}\bigr).\nonumber
\end{eqnarray}
$T_{11}$ can be dealt with in the same way as in the proof of Theorem
\ref{Th1}. It is asymptotically normal with mean $\mu_{1}$ and variance
$\sigma_{1}^{2}$ given in Theorem \ref{Th1}. By the definition, $T_{12}%
=d_{1n}$. We now study the third term $T_{13}$. By (\ref{ap6a}) and
(\ref{A.17b}), $T_{13}$ admits the following decomposition:
\begin{eqnarray*}
\frac{1}{2}T_{13}  &=& \sum_{i=1}^{n}g_{n}(X_{i},Z_{i})w_{i}[A_{n}^{\ast
}(X_{i},Z_{i})-L_{n2}(X_{i},Z_{i})-L_{n1}(X_{i},Z_{i})]
\\
&=& \sum_{i=1}^{n}g_{n}(X_{i},Z_{i})w_{i}\frac{1}{n} \sum_{j=1}^{n}%
\{W_{n}(X_{j}-X_{i},X_{i};h_{1})\varepsilon_{j}^{\ast}(Z_{i})
\\
&&\qquad\hspace*{80pt} {}  -W_{n}(Y_{j}-X_{i};X_{i};h_{3})e_{j+1}(Z_{i})
\\
&&\qquad\hspace*{80pt}\hspace*{43pt}{} -q^{\ast}(X_{i},Z_{j}
)\varepsilon_{j+1}(Z_{i})\}
\\
&&{}  +o_{p}\bigl(1/\sqrt{h_{1}h_{2}}\bigr)+O\bigl(n\delta_{n}(b_{1}^{2}+b_{2}^{2}
)\bigr)+O(\delta_{n}h_{1}^{-1}h_{2}^{-1})
\\
&=& \sum_{i\neq j}\frac{1}{n}g_{n}(X_{i},Z_{i})w_{i}\pi^{-1}(X_{i}
)\{W_{h_{1}}(X_{j}-X_{i})\varepsilon_{j}^{\ast}(Z_{i})
\\
&&\qquad\hspace*{102pt} {} -W_{h_{3}}(Y_{j}-X_{i})e_{j+1}(Z_{i})
\\
&&\qquad\hspace*{102pt}\hspace*{14pt}{} -q^{\ast}(X_{i},Z_{j})\varepsilon
_{j+1}(Z_{i})\}
\\
&&{} +o_{p}\bigl(1/\sqrt{h_{1}h_{2}}\bigr)+O\bigl(n\delta_{n}(b_{1}^{2}+b_{2}^{2}
)\bigr)+O(\delta_{n}h_{1}^{-1}h_{2}^{-1})
\\
&\equiv& \sum_{i\neq j}\varphi(i,j)+o_{p}\bigl(1/\sqrt{h_{1}h_{2}}\bigr) +O\bigl(n\delta
_{n}(b_{1}^{2}+b_{2}^{2})\bigr)+O\bigl(\delta_{n}/(h_{1}h_{2})\bigr).
\end{eqnarray*}
The first term above is a $U$-statistic with the typical element $\varphi
(i,j)$. Let $\varphi^{\ast}(i,j)=\varphi(i,j)+\varphi(j,i)$. Then
$\varphi^{\ast}(i,j)$ is a symmetric kernel and
\[
T_{13}=\sum_{1\leq i<j\leq n}\varphi^{\ast}(i,j)+O\bigl(\delta_{n}/(h_{1}%
h_{2})\bigr)+o_{p}\bigl(1/\sqrt{h_{1}h_{2}}\bigr).
\]
Put $\tilde{\varphi}(i)=\int\varphi^{\ast}(i,j) \, dF_{j}$ and $\tilde{\varphi
}(i,j)=\varphi^{\ast}(i,j)-\tilde{\varphi}(i)-\tilde{\varphi}(j).$ Then by the
Hoeffding decomposition, we have
\[
\sum_{1\leq i<j\leq n}\varphi^{\ast}(i,j)=\sum_{1\leq i<j\leq n}
\tilde{\varphi}(i,j)+(n-1)\sum_{i=1}^{n}\tilde{\varphi}(i).
\]
It is easy to show that $E[h_{1}h_{2}\tilde{\varphi}(i,j)]^{2(1+\delta)}
=O(\delta_{n}^{2(1+\delta)}n^{-2(1+\delta)}h_{1}h_{2}).$ Therefore, applying
the lemma with $\delta=1$ of \cite{gaoking04},
we obtain that
\[
E\biggl\{\sum_{1\leq i<j\leq n}\tilde{\varphi}(i,j)\biggr\}^{2}=o\bigl(1/(h_{1}%
h_{2})\bigr).
\]
Therefore,
\begin{equation}\label{Th3:eq5a}
T_{13}=(n-1)\sum_{i=1}^{n}\tilde{\varphi}(i)+o_{p}\bigl(1/\sqrt{h_{1}h_{2}}\bigr)
+O\bigl(\delta_{n}/(h_{1}h_{2})\bigr).
\end{equation}
By the definition of $\tilde{\varphi}_{i}$, it can be written that
\begin{eqnarray*}
\tilde{\varphi}(i)  &=& \frac{2}{n}g_{n}(X_{i},Z_{i})w(X_{i},Z_{i})\pi
^{-1}(X_{i})
\\
&&{}\times \int\{W_{h_{1}}(x_{j}-X_{i})\varepsilon_{j}^{\ast}(Z_{i})
  -W_{h_{3}}(y_{j}-X_{i})e_{j+1}(Z_{i})
  \\
  &&\hspace*{134pt}{}-q^{\ast}(X_{i},z_{j})\varepsilon
_{j+1}(Z_{i})\}\, dF_{j}
\\
& \equiv&\tilde{\varphi}_{1}(i)+\tilde{\varphi}_{2}(i)+\tilde{\varphi}_{3}(i),
\end{eqnarray*}
where
\begin{eqnarray*}
\tilde{\varphi}_{1}(i)&=&\frac{2}{n}g_{n}(X_{i},Z_{i})w(X_{i},Z_{i})\pi
^{-1}(X_{i})\int W_{h_{1}}(x_{j}-X_{i})\varepsilon_{j}^{\ast}(Z_{i}) \,dF_{j},
\\
\tilde{\varphi}_{2}(i)&=&-\frac{2}{n}g_{n}(X_{i},Z_{i})w(X_{i},Z_{i})\pi
^{-1}(X_{i})\int W_{h_{3}}(y_{j}-X_{i})e_{j+1}(Z_{i})\, dF_{j}
\end{eqnarray*}
and $\tilde{\varphi}_{3}(i)=-\frac{2}{n}g_{n}(X_{i},Z_{i})w(X_{i},Z_{i}%
)\pi^{-1}(X_{i})\int q^{\ast}(X_{i},z_{j})\varepsilon_{j+1}(Z_{i})\, dF_{j}.$
Then by the Fubini theorem and by taking iterative expectation, $E[\tilde
{\varphi}(i)]=0$. Using the central limit theorem for the $\beta$-mixing
process, we get
\[
\frac{(n-1)}{2\sigma_{1A}}\sum_{i=1}^{n}\tilde{\varphi}(i)\stackrel
{{\mathcal{D}}}{\longrightarrow}{\mathcal{N}}(0,1),
\]
where $\sigma_{1A}^{2}=\frac{1}{4}nE[(n-1)^{2}\tilde{\varphi}^{2}(i)]$. By
directly calculating the integration, it can be shown that
\[
\tilde{\varphi}_{1}(i)=\frac{2}{n}g_{n}(X_{i},Z_{i})w(X_{i},Z_{i}
)[p(Z_{i}|X_{i},2\Delta)-p^{2}(Z_{i}|X_{i},2\Delta)]\bigl(1+o(1)\bigr),
\]
$\tilde{\varphi}_{2}(i)=o(g_{n}(X_{i},Z_{i})/n)$ and $\tilde{\varphi}_{3}
(i)=o(g_{n}(X_{i},Z_{i})/n)$. Therefore,
\begin{eqnarray*}
\sigma_{1A}^{2}&=&nE[g_{n}^{2}(X_{1},Z_{1})w^{2}(X_{1},Z_{1})\{p(Z_{i}%
|X_{i},2\Delta)-p^{2}(Z_{i}|X_{i},2\Delta)\}^{2}]
\\
&&{}+o\bigl(1/(h_{1}h_{2})\bigr).
\end{eqnarray*}
By straightforward calculation, it can be shown that the covariance between
$T_{11}$ and $T_{13}$ can be ignored. It follows that the result of the
theorem holds.
\end{pf*}

\begin{pf*}{Proof of Theorem \ref{Th4}} (i) For any given small $\eta>0$, when
$d$ is small enough, $|d_{1n}/\sigma_{1n}|\le\eta$ and $\sigma_{1n}=\sigma
_{1}(1+o(1)).$ Under $H_{0}$, with the selected bandwidths,
\begin{eqnarray*}
(T_{1}-\mu_{1})/\sigma_{1}=O_{p}(1).
\end{eqnarray*}
Therefore, the sequence of critical values $c_{\alpha}$ (depending on $n$) is
bounded in probability. Similarly, under $H_{1n}$, with the selected
bandwidths,
\begin{eqnarray}\label{Th4:eq2}
(T_{1}-\mu_{1}-d_{1n})/\sigma_{1n}=O_{p}(1).
\end{eqnarray}
Note that
\begin{eqnarray*}
P\{(T_{1}-\mu_{1})/\sigma_{1}>c_{\alpha}|H_{1n}\}  &=&P\{(T_{1}-\mu
_{1}-d_{1n})/\sigma_{1n}>(c_{\alpha}\sigma_{1}-d_{1n})/\sigma_{1n}|H_{1n}\}
\\
&\le& P\{(T_{1}-\mu_{1}-d_{1n})/\sigma_{1n}>c_{\alpha}\sigma_{1}/\sigma
_{1n}-\eta|H_{1n}\}.
\end{eqnarray*}
It follows from Theorem \ref{Th3} and Slutsky's theorem that
\begin{eqnarray*}
\limsup_{d\to0}\limsup_{n\to\infty} P\{(T_{1}-\mu_{1})/\sigma_{1}\ge
c_{\alpha}|H_{1n}\}\le\alpha.
\end{eqnarray*}

(ii) For any given $M>0$, by taking $d$ sufficiently large, there exists an
$N$, when $n>N$, $d_{1n}/\sigma_{1n}\ge M.$ Therefore,
\begin{eqnarray*}
P\{(T_{1}-\mu_{1})/\sigma_{1}>c_{\alpha}|H_{1n}\}  &  \ge P\{(T_{1}-\mu
_{1}-d_{1n})/\sigma_{1n}>c_{\alpha}\sigma_{1}/\sigma_{1n}-M|H_{1n}\}.
\end{eqnarray*}
By (\ref{Th4:eq2}), we have
\begin{eqnarray*}
\liminf_{d\to\infty}\liminf_{n\to\infty} P\{(T_{1}-\mu_{1})/\sigma
_{1}>c_{\alpha}|H_{1n}\}=1.
\end{eqnarray*}
\upqed\end{pf*}

\begin{pf*}{Proof of Theorems \ref{Th5} and \ref{Th6}} We put the proofs in the
supplemental materials \cite{2a}.
\end{pf*}

\section*{Acknowledgments} The authors thank the Associate Editor and the referees for
constructive comments that substantially improved an earlier version
of this paper.

%

\begin{supplement}
\sname{Supplement}\label{suppA}
\stitle{Additional technical details}
\slink[doi]{10.1214/09-AOS763SUPP}
\sdatatype{.pdf}
\slink[url]{}
\sdescription{We provide detailed proofs for Lemmas \ref{le0}--\ref{lead3} and Theorems \ref{Th5}--\ref{Th6}.
Modern nonparametric smoothing techniques and theory of $U$-statistics are used.}
\end{supplement}

\printaddresses


\begin{thebibliography}{99}

\bibitem{yacrfs96}
\textsc{A{\"{\i}}t-Sahalia, Y.} (1996).
Testing continuous-time models of the spot interest rate.
\emph{Review of Financial Studies} \textbf{9} 385--426.


\bibitem{2a}
\textsc{A{\"{\i}}t-Sahalia, Y., Fan, J.} and \textsc{Jiang, J.} (2010).
Supplement to ``Nonparametric tests of the Markov hypothesis in continuous-time models.''
\href{http://dx.doi.org/10.1214/09-AOS763SUPP}{DOI: 10.1214/09-AOS763SUPP}.

\bibitem{yacfanpeng09}
\textsc{A{\"{\i}}t-Sahalia, Y.}, \textsc{Fan, J.} and  \textsc{Peng,
  H.} (2009).
Nonparametric transition-based tests for jump-diffusions.
\textit{J. Amer. Statist. Assoc.} \textbf{104}
  1102--1116.

\bibitem{azzalinibowmanhardle89}
\textsc{Azzalini, A.}, \textsc{Bowman, A.~N.} and  \textsc{H\"{a}rdle,
  W.} (1989).
On the use of nonparametric regression for model checking.
\textit{Biometrika} \textbf{76} 1--11.
\MR{0991417}

\bibitem{bickelrosenblatt93}
\textsc{Bickel, P.~J.} and  \textsc{Rosenblatt, M.} (1973).
On some global measures of the deviation of density function
  estimates.
\textit{Ann. Statist.} \textbf{1} 1071--1095.
\MR{0348906}

\bibitem{caverhill94}
\textsc{Caverhill, A.} (1994).
When is the short rate {M}arkovian?
\textit{Math. Finance} \textbf{4} 305--312.
\MR{1299241}

\bibitem{chengao04b}
\textsc{Chen, S.~X.} and  \textsc{Gao, J.} (2004).
On the use of the kernel method for specification tests of diffusion
  models.
Technical report, Iowa State Univ.

\bibitem{chengaotong08}
\textsc{Chen, S.~X.}, \textsc{Gao, J.} and  \textsc{Tang, C.}
  (2008).
A test for model specification of diffusion processes.
\textit{Ann. Statist.} \textbf{36} 167--198.
\MR{2387968}

\bibitem{cir85}
\textsc{Cox, J.~C.}, \textsc{Ingersoll, J.~E.} and  \textsc{Ross,
  S.~A.} (1985).
A theory of the term structure of interest rates.
\textit{Econometrica} \textbf{53} 385--408.
\MR{0785475}

\bibitem{fan96}
\textsc{Fan, J.} (1996).
Test of significance based on wavelet thresholding and {N}eyman's
  truncation.
\textit{J. Amer. Statist. Assoc.} \textbf{91}
  674--688.
\MR{1395735}

\bibitem{fanjiang05}
\textsc{Fan, J.} and  \textsc{Jiang, J.} (2005).
Generalized likelihood ratio tests for additive models.
\textit{J. Amer. Statist. Assoc.} \textbf{100}
  890--907.
\MR{2201017}

\bibitem{fanjiang07}
\textsc{Fan, J.} and  \textsc{Jiang, J.} (2007).
Nonparametric inference with generalized likelihood ratio tests (with
  discussion).
\textit{Test} \textbf{16} 409--478.
\MR{2365172}

\bibitem{fanyao03}
\textsc{Fan, J.} and  \textsc{Yao, Q.} (2003).
\textit{Nonlinear Time Series: Nonparametric and Parametric Methods}.
Springer, New York.
\MR{1964455}

\bibitem{fanyaotong96}
\textsc{Fan, J.}, \textsc{Yao, Q.} and  \textsc{Tong, H.} (1996).
Estimation of conditional densities and sensitivity measures in
  nonlinear dynamical systems.
\textit{Biometrika} \textbf{83} 189--206.
\MR{1399164}

\bibitem{fanyim04}
\textsc{Fan, J.} and  \textsc{Yim, T.-H.} (2004).
A data-driven method for estimating conditional densities.
\textit{Biometrika} \textbf{91} 819--834.
\MR{2126035}

\bibitem{fanzhangzhang01}
\textsc{Fan, J.}, \textsc{Zhang, C.} and  \textsc{Zhang, J.} (2001).
Generalized likelihood ratio statistics and {W}ilks phenomenon.
\textit{Ann. Statist.} \textbf{29} 153--193.
\MR{1833962}

\bibitem{fouquepapanicolaousircar00}
\textsc{Fouque, J.-P.}, \textsc{Papanicolaou, G.} and  \textsc{Sircar,
  K.~R.} (2000).
\textit{Derivatives in Financial Markets with Stochastic Volatility}.
Cambridge Univ. Press, London.
\MR{1768877}

\bibitem{gaocasas08}
\textsc{Gao, J.} and  \textsc{Casas, I.} (2008).
Specification testing in discretized diffusion models: {T}heory and
  practice.
\textit{J. Econometrics} \textbf{147} 131--140.
\MR{2472987}

\bibitem{gaoking04}
\textsc{Gao, J.} and  \textsc{King, M.} (2004).
Model specification testing in nonparametric and semiparametric time
  series econometrics.
Technical report,  Univ. Western Australia.

\bibitem{hallracineli04}
\textsc{Hall, P.}, \textsc{Racine, J.} and  \textsc{Li, Q.} (2004).
Cross-validation and the estimation of conditional probability
  densities.
\textit{J. Amer. Statist. Assoc.} \textbf{99}
  1015--1026.
\MR{2109491}

\bibitem{hardlemammen93}
\textsc{H{\"{a}}rdle, W.} and  \textsc{Mammen, E.} (1993).
Comparing nonparametric versus parametric regression fits.
\textit{Ann. Statist.} \textbf{21} 1926--1947.
\MR{1245774}

\bibitem{heathjarrowmorton92}
\textsc{Heath, D.}, \textsc{Jarrow, R.} and  \textsc{Morton, A.}
  (1992).
Bond pricing and the term structure of interest rates: A new
  methodology for contingent claims evaluation.
\textit{Econometrica} \textbf{60} 77--105.

\bibitem{heston93}
\textsc{Heston, S.} (1993).
A closed-form solution for options with stochastic volatility with
  applications to bonds and currency options.
\textit{Review of Financial Studies} \textbf{6} 327--343.


\bibitem{hongli05}
\textsc{Hong, Y.} and  \textsc{Li, H.} (2005).
Nonparametric specification testing for continuous-time models with
  applications to term structure of interest rates.
\textit{Review of Financial Studies} \textbf{18} 37--84.

\bibitem{hyndmanyao02}
\textsc{Hyndman, R.} and  \textsc{Yao, Q.} (2002).
Nonparametric estimation and symmetry tests for conditional density
  functions.
\textit{J. Nonparametr. Statist.} \textbf{14} 259--278.
\MR{1905751}

\bibitem{ingster93}
\textsc{Ingster, Y.} (1993).
Asymptotically minimax hypothesis testing for nonparametric
  alternatives {I--III}.
\textit{Math. Methods  Statist.} \textbf{2} 85--114;
 \textbf{3} 171--189; \textbf{4} 249--268.


\bibitem{lepskispokoiny99}
\textsc{Lepski, O.} and  \textsc{Spokoiny, V.} (1999).
Minimax nonparametric hypothesis testing: {T}he case of an
  inhomogeneous alternative.
\textit{Bernoulli} \textbf{5} 333--358.
\MR{1681702}


\bibitem{revuzyor}
\textsc{Revuz, D.} and  \textsc{Yor, M.} (1994).
\textit{Continuous Martingales and Brownian Motion}, 2nd ed.
Springer, Berlin.
\MR{1303781}

\bibitem{spokoiny96}
\textsc{Spokoiny, V.~G.} (1996).
Adaptive hypothesis testing using wavelets.
\textit{Ann. Statist.} \textbf{24} 2477--2498.
\MR{1425962}

\bibitem{vasicek77}
\textsc{Vasicek, O.} (1977).
An equilibrium characterization of the term structure.
\textit{Journal of Financial Economics} \textbf{5} 177--188.

\end{thebibliography}
\end{document}